\documentclass[microtype]{gtpart} 
\usepackage{geometry}
\usepackage[parfill]{parskip}    
\usepackage{amsfonts} 
\usepackage{amsmath} 
\usepackage{amsthm}
\usepackage{dsfont}
\usepackage{epstopdf}
\usepackage{hyperref} 
\usepackage[all]{xy} 
\usepackage[pdftex]{color,graphicx}
\usepackage{psfrag}

\title{\ \ \ \ \ \ \ \ \ Heegaard Floer homology of broken fibrations over the circle}
\author{Yank\i\ Lekili}
\givenname{Yank\i} \surname{Lekili}
\address{Department of Mathematics, M.I.T., Cambridge MA 02139, USA}
\email{lekili@math.mit.edu}


\begin{abstract} We extend Perutz's Lagrangian matching invariants to
3--manifolds which are not necessarily fibred using the technology of
holomorphic quilts. We prove an isomorphism of these invariants with
Ozsv\'ath-Szab\'o's Heegaard Floer invariants for certain extremal spin$^c$
structures. As applications, we give new calculations of Heegaard Floer
homology of certain classes of 3--manifolds, and a characterization of
Juh\'asz's sutured Floer homology.
\end{abstract}

\newcommand{\f}[1]{\mathbb{#1}} 
\newcommand{\s}{\mathfrak{s}}
\newcommand{\x}{\mathbf{x}}
\newcommand{\y}{\mathbf{y}}
\newcommand{\av}{\mathbf{a}}
\newcommand{\bv}{\mathbf{b}}
\newcommand{\lft}{\text{left}} 
\newcommand{\rgh}{\text{right}}
\newcommand{\ten}{\otimes} 
\newcommand{\del}{\partial}
\newcommand{\delbar}{\bar{\partial}} 
\newcommand{\sym}{\text{Sym}}
\newcommand{\hilb}{\text{Hilb}}
\newcommand{\OO}[1]{\operatorname{O}(#1)}

\newcommand{\QED}{\vspace{-.31in}\begin{flushright}\qed\end{flushright}}

\newtheorem{theorem}{Theorem} 
\newtheorem{lemma}[theorem]{Lemma}
\newtheorem{proposition}[theorem]{Proposition}
\newtheorem{definition}[theorem]{Definition}
\newtheorem{corollary}[theorem]{Corollary} 

\begin{document}

\maketitle
\section{Introduction}

In this paper, we study two seemingly different Floer theoretical invariants of three-
and four-manifolds. These are Perutz's Lagrangian matching invariants and Ozsv\'ath
and Szab\'o's Heegaard Floer theoretical invariants. The main result of this paper is
an isomorphism between the 3--manifold invariants of these theories for certain
spin$^c$ structures, namely {\it quilted Floer homology} and {\it Heegaard Floer
homology}. We also outline how the techniques here can be generalized to obtain an
identification of 4--manifold invariants.

Before giving a review of both of the above mentioned theories, we give the definition
of a broken fibration over $S^1$, which will be an important part of the topological
setting that we will be working with.

\begin{definition} A map $f:Y\to S^1$ from a closed oriented smooth $3$--manifold $Y$
	to $S^1$ is called a broken fibration if $f$ is a circle-valued Morse function
	with all of the critical points having index $1$ or $2$. 
\end{definition}

The terminology is inspired from the terminology of broken Lefschetz fibrations on
$4$--manifolds, to which we will return later in this paper in Section \ref{4manifolds}.
We remark that a $3$--manifold admits a broken fibration if and only if $b_1(Y)>0$,
and if it admits one, it admits a broken fibration with connected fibres. 

We will restrict ourselves to broken fibrations with connected fibres and we will denote by
$\Sigma_{\text{max}}$ and $\Sigma_{\text{min}}$ two fibres with maximal and minimal
genus respectively. We denote by $\mathcal{S}(Y|\Sigma_{\text{min}})$, the spin$^c$ structures $\s$
on $Y$ such that $\langle c_1(\s) , [\Sigma_{\text{min}}] \rangle =
\chi(\Sigma_{\text{min}})$ (those spin$^c$ structures which satisfy the adjunction
{\it equality} with respect to the fibre with minimal genus).

\begin{definition} \label{novikov}The universal Novikov ring $\Lambda$ over $\f{Z}_2$ is the ring of formal power
	series $\Lambda = \Sigma_{ r\in \f{R}} \ a_r t^r $ with $a_r \in \f{Z}_2$ such
	that $\#\{r | a_r\neq 0, r < N\} < \infty $ for any $N\in \f{R}$. 
\end{definition}	

We first give a definition of a new invariant $QFH'(Y,f,\s;\Lambda)$ for all spin$^c$
structures in $\mathcal{S}(Y|\Sigma_{\min})$ and prove an isomorphism between this
variant of quilted Floer homology of a broken fibration $f:Y\to S^1$ (with
coefficients in the universal Novikov ring) and the Heegaard Floer homology of $Y$
perturbed by a closed 2-form $\eta$ that pairs positively with the fibers of $f$:

\begin{theorem} \label{maintheorem} $\text{QFH}'(Y, f, \s; \Lambda) \simeq
	HF^{\pm}(Y, \eta, \s)$ for $\s \in \mathcal{S}(Y|\Sigma_{\text{min}})$.

\end{theorem}

When $g(\Sigma_\text{min})$ is at least $2$, the coefficients can be taken to be in $\f{Z}_2$ (in this case admissibility of our diagrams are automatic, therefore we
do not need to use perturbations).

\begin{corollary}  
	Suppose that $g(\Sigma_\text{min}) >1$. Then for
	$\s\in\mathcal{S}(Y|\Sigma_\text{min})$ we have \[ QFH'(Y,f,\s;\f{Z}_2)
	\simeq HF^+(Y,\s; \f{Z}_2) \]
\end{corollary}

As corollaries of this result, we give new calculations of Heegaard Floer homology groups for
certain manifolds for which $QFH'(Y,f,\s)$ is easy to calculate. We give several such
calculations among which the following is particularly interesting. 

\begin{corollary}
Suppose $f$ has only two critical points, and let $\alpha,\beta \subset
\Sigma_\text{max}$ be the vanishing cycles of these critical points. Then, $\oplus_{\s
\in \mathcal{S}(Y|\Sigma_\text{min})} HF^+(Y,\eta,\s)$ is free of rank
$\iota(\alpha,\beta)$, the geometric intersection number between $\alpha$ and $\beta$. 
Furthermore, if $g(\Sigma_\text{min})>1$ then the result holds over $\f{Z}_2$, i.e.
\[ \oplus_{\s \in \mathcal{S}(Y|\Sigma_\text{min})} HF^+(Y,\s) = \f{Z}_2^{\iota(\alpha,\beta)} \]
\end{corollary} 
 
The second main theorem proves that the invariants $QFH'(Y,f,\s;\Lambda)$ that we
defined are isomorphic to the quilted Floer homology groups coming from Perutz's theory
of Lagrangian matching invariants. Unlike $QFH'(Y,f,\s;\Lambda)$, for technical reasons these latter invariants are only defined in the
case $g(\Sigma_{\text{max}}) < 2g(\Sigma_{\text{min}})$. Thus, we have the following
theorem :

\begin{theorem} 
\label{maintheorem2} Suppose that $Y$ admits a broken fibration with 
$g(\Sigma_{\text{max}}) < 2g(\Sigma_{\text{min}})$. Then for $\s \in
\mathcal{S}(Y|\Sigma_{\text{min}})$, $QFH(Y,f,\s;\Lambda)$ is
well-defined and
$$\text{QFH}'(Y, f, \s; \Lambda) \simeq QFH(Y, f, \s;\Lambda)$$  \end{theorem}

As before, we have the same result over $\f{Z}_2$ when $g(\Sigma_\text{min})$ is at least $2$. 

In Section \ref{section2}, we construct a Heegaard diagram associated with a broken
fibration and investigate the properties of this diagram. We also give a calculation
of perturbed Heegaard Floer homology of fibred 3--manifolds for
$\s\in\mathcal{S}(Y|F)$. In Section \ref{section3}, we give a definition of quilted
Floer homology in the language of Heegaard Floer theory and prove that it is
isomorphic to the Heegaard Floer homology for the spin$^c$ structures under
consideration. Here, we give several corollaries of our first main result, including new
calculations of Heegaard Floer homology groups and a characterization of Juh\'asz's
sutured Floer homology. In Section \ref{section4}, we give a complete definition of $QFH(Y,f,\s)$ and we prove our second main theorem,
namely that the group defined in Section \ref{section3} is isomorphic to the original
definition of quilted Floer homology in terms of holomorphic quilts. Finally in Section
\ref{4manifolds} we discuss the extension of this isomorphism to four-manifold invariants. 

We now proceed to review the theories and the notation that are involved in our
theorem.  

\subsection{(Perturbed) Heegaard Floer homology} 
\label{introperturb}

In this section, we review the construction of Heegaard Floer homology, introduced by
Ozsv\'ath and Szab\'o \cite{OS}. The usual construction involves certain admissibility conditions, however there is a variant of
Heegaard Floer homology where Novikov rings and perturbations by closed $2$-forms
are introduced in order to make the Heegaard Floer homology group well-defined
without any admissibility condition. Our account will be brief since this theory has
been well developed in the literature. The reader is encouraged to turn to
\cite{perturbed} for a more detailed account of perturbed Heegaard Floer theory.
Furthermore, we will mostly find it convenient to work in the set up of Lipshitz's
cylindrical reformulation of Heegaard Floer homology \cite{lipshitz}. 

Let $(\Sigma_{g}, \boldsymbol{\alpha},\boldsymbol{\beta},z)$ be a pointed Heegaard diagram of a
$3$--manifold $Y$.  This gives rise to a pair of Lagrangian tori $\mathbb{T}_\alpha$,
$\mathbb{T}_\beta$ in $\sym^g(\Sigma_g)$, together with a holomorphic hypersurface $Z=
z \times \sym^{g-1}(\Sigma_g)$. The Heegaard Floer homology of $Y$ is the Lagrangian
Floer homology of these tori, where one uses the orbifold symplectic form pushed down
from $\Sigma^{\times g}_g$, though one can also use honest symplectic forms (see
\cite{handleslide}).
The differential is twisted by keeping track of the intersection number $n_z$ of holomorphic disks contributing to the differential with $Z$.  More precisely, the Heegaard Floer chain complex $CF^+(Y)$ is
freely generated over $\f{Z}$ by $[\x,i]$ where $\x$ is an intersection point of
$\mathbb{T}_\alpha$ and $\mathbb{T}_\beta$ and $i\in\f{Z}_{\geq 0}$, and the
differential is given by \[ \del^+([\x,i]) = \sum_{\y} \sum_{\varphi\in\pi_2(\x,\y),
n_z(\varphi)\leq i } \#\widehat{\mathcal{M}}(\varphi) [y,i-n_z(\varphi)] \]

where $\#\widehat{\mathcal{M}}(\varphi)$, as usual in Lagrangian Floer homology, refers to a count of holomorphic disks with boundary on $\f{T}_\alpha$ and $\f{T}_\beta$ connecting $\x$ and $\y$. The above definition only makes sense under certain admissibility conditions so that
the sum on the right hand side of the differential is finite. In general, one can consider a twisted
version of the above chain complex by a closed $2$-form in $\Omega^2(Y)$. This is
called the {\it perturbed Heegaard Floer homology}. The chain complex $CF^+(Y,\eta)$
is freely generated over $\Lambda$ (see Definition \ref{novikov}) by $[\x,i]$ where
$\x$ is an intersection point and $i$ is a nonnegative integer as before, and the
differential is twisted by the area $\int_{[\varphi]} \eta$ of the holomorphic disks
that contribute to the differential. More precisely, the differential of the perturbed
theory is given by \[ \del^+([\x,i]) = \sum_{\y}
   \sum_{\varphi\in\pi_2(\x,\y), n_z(\varphi)\leq i }
   \#\widehat{\mathcal{M}}(\varphi) t^{\eta(\varphi)} [y,i-n_z(\varphi)] \]

Note that if $\varphi_1,\varphi_2$ are two holomorphic discs that connect an
intersection point $\x$ to $\y$, then their difference is a periodic domain $P$ and we
have the equality $\eta(\varphi_1) - \eta(\varphi_2) = \eta ([P])$, where the latter
only depends on the cohomology class of $\eta$. We remark that although the
differential depends on the choice of a representative of the class $[\eta]$, the
isomorphism class of the homology groups is determined by $\text{Ker}(\eta) \cap
H_2(Y;\f{Z})$. 

Recall that a $2$--form is said to be generic when $\text{Ker}(\eta) \cap H_2(Y;\f{Z})
=\{ 0 \}$.  For a generic form coming form an area form on the Heegaard surface, $HF^+(Y,\eta)$ is defined without any admissibility conditions on the Heegaard diagram.

\subsection{Quilted Floer homology of a $3$--manifold}

In this section, we review the definition of quilted Floer homology of a $3$--manifold
$Y$ equipped with a broken fibration $f:Y\to S^1$. The general theory of holomorphic
quilts is under systematically developed by Wehrheim and Woodward \cite{katrin}, though
the case we consider also appears in the work of Perutz \cite{gysin}. The relevant part of the
theory in the setting of $3$-manifolds is obtained from Perutz's construction of
Lagrangian matching conditions associated with critical values of broken fibrations,
which we now review from \cite{LM1}. 

Given a Riemann surface $(\Sigma, j)$ and an embedded circle $L \subset \Sigma $,
denote by $\Sigma_L$ the surface obtained from $\Sigma$ by surgery along $L$, i.e., by
removing a tubular neighborhood of $L$ and gluing in a pair of discs. To such data,
Perutz associates a distinguished Hamiltonian isotopy-class of Lagrangian
correspondences ${V}_L \subset \text{Sym}^n(\Sigma) \times \text{Sym}^{n-1}(\Sigma_L)$
(where the symmetric products are equipped with K\"ahler forms in suitable cohomology
classes, see \cite{LM1}). These are described in terms of a symplectic degeneration of
$\sym^n(\Sigma)$. More precisely, one considers an elementary Lefschetz fibration over
$D^2$ with regular fibre $\Sigma$ and a unique vanishing cycle $L$ which collapses at
the origin. Then one passes to the relative Hilbert scheme, $\hilb^n_{D^2}(\Sigma)$,
of this fibration (the resolution of the singular variety obtained by taking
fibre-wise symmetric products). The regular fibres of the induced map from
$\hilb^n_{D^2}(\Sigma)$ are identified with $\sym^n(\Sigma)$, and the fibre above the
origin has a codimension $2$ singular locus which can be identified with
$\sym^{n-1}(\Sigma_L)$. $V_{L}$ then arises as the vanishing cycle of this fibration.

Given a $3$--manifold $Y$ and a broken fibration $f:Y\to S^1$, the {\it quilted Floer
homology} of $Y$, $QFH(Y,f)$, is a Lagrangian intersection theory graded by
$\text{spin}^c$ structures on $Y$. Let $p_1, p_2,\ldots, p_k$ be the set of critical
values of $f$. Pick points $p^{\pm}_i$ in a small neighborhood of each $p_i$ so that
the fibre genus increases from $p^-_i$ to $p^+_i$. For $\s \in \text{spin}^c (Y)$, let
$\nu : S^1 \backslash \text{crit}(f) \to \f{Z}_{\geq 0} $ be the locally constant
function defined by $ \langle c_1(\s), [F_s] \rangle = 2\nu(s) + \chi(F_s) $, where
$F_s = f^{-1}(s)$. Then the construction in the previous paragraph gives Lagrangian
correspondences $L_{p_i} \subset \sym^{\nu(p^+_i)}(F_{p^+_i}) \times
\sym^{\nu(p^-_i)}(F_{p^-_i})$. The quilted Floer homology of $Y$,
$QFH(L_{p_1},\ldots,L_{p_k})$, is then generated by horizontal (with respect to the
gradient flow of $f$) multi-sections of $f$ which match along the Lagrangians
$L_{p_1},\ldots,L_{p_k}$ at the critical values of $f$, and the differential counts
rigid holomorphic ``quilted cylinders'' connecting the generators, \cite{gysin},
\cite{katrin} (see Section \ref{qfhdef} for a detailed definition).

There are various technical difficulties involved in the definition of $QFH(Y,f,\s)$
due to bubbling of holomorphic curves. These are addressed by different means
depending on the value of $\langle c_1(\s), [\Sigma_{\text{max}}] \rangle$. The easiest
case is the (positively) monotone case, that is when  $\langle c_1(\s), [\Sigma_{\text{max}}] \rangle
> 0$, where holomorphic bubbles are a priori excluded.  However, for
$\s\in\mathcal{S}(Y|\Sigma_{\text{min}})$ we will almost never be in the monotone
case. In the strongly negative case, that is when $\langle c_1(\s),
[\Sigma_{\text{max}}] \rangle \leq \chi(\Sigma_{\text{max}})/2$, one can still eliminate
bubbles a priori by standard means. For the rest of the cases, bubbles might and will
occur in general, therefore complications arise. One then tries to
establish a proper combinatorial rule for handling bubbled configurations.  One could
also try to use the more technical machinery of \cite{liuTian} or \cite{cieMohn} in order
to tackle this case. Another related issue is showing that quilted Floer homology is
an invariant of a three manifold. The isomorphism constructed in this paper shows this
in an indirect way for the spin$^c$ structures under consideration. We will return to
this question and various well-definedness questions in \cite{lekiliPerutz}. 

In this paper, we will deal with the spin$^c$ structures $\s \in
\mathcal{S}(Y|\Sigma_{\text{min}})$. In this case, quilted Floer homology has been
defined only in the strongly negative case, which is equivalent to requiring
$g(\Sigma_\text{max}) < 2 g(\Sigma_{\text{min}})$ (see Section \ref{section4} for details). However, we will define a variant
of quilted Floer homology, which we will denote by $QFH'(Y,f,\s)$ that suits our
purposes and avoids these technical issues, hence is well-defined in all cases; see
Section \ref{variant} for the definition. We will prove that in the case when
$QFH(Y,f,\s)$ is defined, it is isomorphic to $QFH'(Y,f,\s)$ (this is the content of
our Theorem \ref{maintheorem2} above). Then, Theorem \ref{maintheorem}, which
establishes an isomorphism between $QFH'(Y,f,\s)$ and $HF^+(Y,\s)$, will show that
$QFH(Y,f,\s)$, when defined, is isomorphic to $HF^+(Y,\s)$.

Finally, we remark that in the case when $f:Y\to S^1$ is a fibration, $QFH(Y,f)$ is
given as a fixed point Floer homology theory on the moduli space of vortices and was
first introduced by Salamon in \cite{salamon}. In this case, the spin$^c$ structures
$\s\in\mathcal{S}(Y|\Sigma)$ corresponds to taking the zeroth symmetric product of the
fibres. In this case, it is natural to set $QFH(Y, f) = \Lambda$ if $\s$ is the canonical tangent spin$^c$ structure, and $QFH(Y,f)= 0$ for other $\s \in \mathcal{S}(Y|\Sigma)$. 
\paragraph{Acknowledgements}

This paper is largely a rewrite of part of the author's PhD thesis. The author
would like to thank his advisor Denis Auroux for his generosity with time and
ideas throughout three years. Special thanks to Tim Perutz for explaining the
details of his work on which this paper builds on, and Robert Lipshitz for a
critical discussion. He is also indebted to Matthew Hedden, Max Lipyanskiy,
Peter Ozsv\'ath for helpful discussions.  Thanks to Peter Kronheimer and Tomasz
Mrowka for their interest in this work. Finally, thanks to the referee for
their helpful comments and suggestions. This work was partially funded by NSF
grants DMS-0600148 and DMS-0706967.

\section{Heegaard diagram for a broken fibration on Y}
\label{section2}
\subsection{A standard Heegaard diagram}
\label{section2.1}
We start with a $3$--manifold $Y$ with $b_1 > 0$. Then $Y$ admits a broken fibration over
$S^1$.  Consider such a Morse function $f: Y\to S^1$ with the following additional
properties : 

\begin{itemize}
	\item $F_{-1}=\Sigma_{\text{max}}$ has the maximal genus $g_{\text{max}} = g$
		and $F_1 = \Sigma_{\text{min}}$ has the minimal genus
		$g_{\text{min}}=k$ among fibres of $f$. 
	\item The fibres are connected.
	\item The genera of the fibres are in decreasing order as one travels
		clockwise and counter-clockwise from $-1$ to $1$. 
\end{itemize}

It is easy to see that a broken fibration with these properties exists if and
only if $b_1>0$. In fact, any broken fibration with connected fibers can be
deformed into one with these properties by an isotopy that changes the order of
the critical values. 

\begin{figure}[!h]
	\centering
	\includegraphics{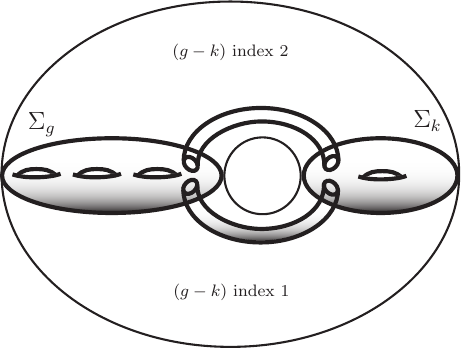}  \caption{Heegaard surface for a broken
	fibration}
\label{heegaardSurface}
\end{figure}

We will now construct a Heegaard diagram for $Y$ adapted to $f$. Roughly speaking, the
Heegaard surface $\Sigma$ will be obtained by connecting $\Sigma_{\text{max}}$ and
$\Sigma_{\text{min}}$ by two ``tubes'' traveling clockwise and counter-clockwise from
$\Sigma_{\text{max}}$ to $\Sigma_{\text{min}}$.  More precisely, start with a section
$\gamma$ of $f$ over $S^1$. Then we can pick a metric for which $\gamma$ is a gradient
flow line of $f$, and since $\gamma$ is disjoint from the critical points of $f$, it
also avoids the stable/unstable manifolds of the critical points. Now pick two
distinct points
$p$ and $q$ on $\Sigma_{\text{max}}$ sufficiently close to the point where $\gamma$
intersects $\Sigma_{\text{max}}$, connect $p$ to $\Sigma_{\text{min}}$ by the
gradient flow line above the northern semi-circle in the base $S^1$ which connects $-1$ to $1$ in the
clockwise direction and connect $q$ to $\Sigma_{\text{min}}$ by the gradient flow line
above the southern
semi-circle, avoiding the critical points of $f$ in both cases. Denote these flow
lines by
$\gamma_p$ and $\gamma_q$ and their end points in $\Sigma_{\text{min}}$ by $\bar{p}$
and $\bar{q}$. Then the Heegaard surface that we are interested in is obtained by
removing discs around $p$, $q$, $\bar{p}$ and $\bar{q}$ and connecting
$\Sigma_{\text{max}}$ to $\Sigma_{\text{min}}$ along $\gamma_p$ and $\gamma_q$ (see
Figure \ref{heegaardSurface}). We denote the resulting surface by \[ \Sigma =
\Sigma_{\text{max}} \cup_{\del N(\gamma_p) \cup \del N(\gamma_q)} \Sigma_{\text{min}}
\]

where $N(\gamma_p)$ and $N(\gamma_q)$ stands for normal neighborhoods of $\gamma_p$ and
$\gamma_q$.

Note that $g(\Sigma) = g+ k + 1$. Denote the point where $\gamma$ intersects
$\Sigma_{\text{max}}$ by $w$ and the point where $\gamma$ intersects
$\Sigma_{\text{min}}$ by $z$. Next, we will describe $\alpha$ and $\beta$ curves
on $\Sigma$ in order to get a Heegaard decomposition of $Y$. First, set $\alpha_0$ to
be $\del N(\gamma_p) \cap f^{-1}(-i)$ and set $\beta_0$ to be $\del N(\gamma_p) \cap
f^{-1}(i)$. The preimage of the
northern semi-circle is a cobordism from $\Sigma_{\text{max}}$ to
$\Sigma_{\text{min}}$ which can be realized by attaching $(g-k)$ 2-handles to
$\Sigma_{\text{max}}\times I$, and hence can be described by the data of $g-k$ disjoint attaching circles on $\Sigma_{\text{max}}$. These we declare to be
$\alpha_1,\ldots,\alpha_{g-k}$.  Similarly the preimage of the southern semi-circle is
a cobordism from $\Sigma_{\text{max}}$ to $\Sigma_{\text{min}}$, encoded by $g-k$ disjoint attaching circles $\beta_1,\ldots,\beta_{g-k}$ on $\Sigma_{\text{max}}$.
Alternatively, these two sets correspond to the stable and unstable manifolds of the
critical points of $f$. More precisely, orienting the base $S^1$ in the clockwise
direction, $\alpha_1,\ldots,\alpha_{g-k}$ are the intersections of the stable
manifolds of the critical points above the northern semi-circle with
$\Sigma_{\text{max}}$, similarly $\beta_1,\ldots,\beta_{g-k}$ are the intersections of
the unstable manifolds of the critical points above the southern semi-circle with
$\Sigma_{\text{max}}$. Note that by choosing $p$ and $q$ sufficiently close to $w$ we
can ensure that they lie in the connected component as $w$ in the complement of $\alpha_1,\ldots,\alpha_{g-k}$ and $\beta_1,\ldots,\beta_{g-k}$. 

\begin{figure}[!h]  
	\centering
	\includegraphics[scale=0.8]{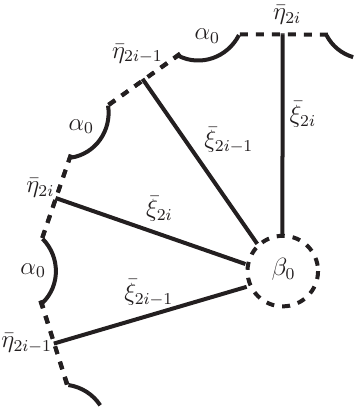}
	\caption{The curves
	$(\bar{\xi}_{2i-1},\bar{\xi}_{2i})$, $(\bar{\eta}_{2i-1},\bar{\eta}_{2i})$}
	\label{curves}
\end{figure}

Next, we describe the remaining curves, $(\alpha_{g-k+1},\ldots, \alpha_{g+k},
\beta_{g-k+1},\ldots, \beta_{g+k})$. Let $F$ be the part of $\Sigma$ which consists of $\Sigma_{\text{max}}$ (except the two discs removed around $p$ and $q$) together with
halves of the connecting tubes up to $\alpha_0$ and $\beta_0$. Thus $F$ is a genus $g$
surface with $2$ boundary components $\alpha_0$ and $\beta_0$. Also, denote by
$\bar{F}$ the complement of $\text{Int}(F)$ in $\Sigma$. Thus $\bar{F}$ is a genus $k$
surface with boundary consisting of $\alpha_0$ and $\beta_0$ and $\Sigma = F
\cup_{\alpha_0 \cup \beta_0} \bar{F}$. Let us also pick $p^+$ and $q^+$ on the
boundary of the disks deleted around $p$ and $q$, and $\bar{p}^+$ and $\bar{q}^+$
their images under the gradient flow (so that they lie on the boundary of the
discs deleted around $\bar{p}$ and $\bar{q}$). Now we can find two $2k$-tuples of ``standard'' pairwise disjoint arcs in $\bar{F}$, $(\bar{\xi}_1,\dots,\bar{\xi}_{2k})$,
$(\bar{\eta}_1,\ldots,\bar{\eta}_{2k})$ such that $\bar{\xi}_i$ intersect
$\bar{\eta}_j$ only if $i=j$, in which case the intersection is transverse at one
point.  Furthermore, we can arrange that the points $z$, $\bar{p}^+$ and $\bar{q}^+$ lie in the same connected component in the complement of these arcs in $\bar{F}$. A nice visualization of these curves on $\bar{F}$ can be obtained by considering a representation of $\bar{F}$ by a $4k$-sided polygon. First, represent a genus $k$ surface by gluing the sides of $4k$-gon in the
way prescribed by the labeling $a_1 b_1 a_1^{-1} b_1^{-1} \ldots a_k b_k a_k^{-1}
b_k^{-1}$ of the sides starting from a vertex and labeling in the clockwise direction.
Now remove a neighborhood of each vertex of the polygon and a neighborhood of a point in its interior. This now represents a genus $k$ surface with two boundary
components. Let us put $\beta_0$ at the boundary of the interior puncture and
$\alpha_0$ at the boundary near the vertices then the curves
$(\bar{\xi}_{2i-1},\bar{\xi}_{2i})$ coincide with the portions of the edges labelled $(a_i, b_i)$ left after removing a neighborhood of each vertex and the
curves $(\bar{\eta}_{2i-1},\bar{\eta}_{2i})$ connect the midpoints of
$(\bar{\xi}_{2i-1},\bar{\xi}_{2i})$ radially to $\beta_0$, see Figure \ref{curves}.

Now, using the gradient flow of $f$ we can flow the arcs
$(\bar{\xi}_1,\ldots,\bar{\xi}_{2k})$ above the northern semi-circle to obtain
disjoint arcs $(\xi_1,\ldots,\xi_{2k})$ in $F$ which do not intersect with
$\alpha_1,\ldots,\alpha_{g-k}$. (Generic choices ensure that the gradient flow does
not hit any critical points.) The flow sweeps out discs in $Y$ which bound
$(\alpha_{g-k+1},\ldots,\alpha_{g+k}) = (\xi_1 \cup \bar{\xi}_1,\ldots, \xi_{2k} \cup
\bar{\xi}_{2k})$. Similarly, we define $(\beta_{g-k+1},\ldots,\beta_{g+k})$ by flowing
the arcs $(\bar{\eta}_1,\ldots,\bar{\eta}_{2k})$ above the southern semi-circle. To
complete the Heegaard decomposition of $(Y,f)$ we set the base point on $\bar{F}$ to
be $z$ which lies in the same region as $\bar{p}^+$ and $\bar{q}^+$. Therefore, we
constructed a Heegaard decomposition of $(Y,f)$. We will also make use of a filtration
associated with the base point $w$ which we can ensure to be located in the same
region as $p^+$ and $q^+$ by picking $p$ and $q$ sufficiently close to $w$, which is
the image of $z$ under the gradient flow above the northern and southern semi-circles.
Roughly speaking, this point will be used to keep track of the domains passing through
the connecting ``tubes''.

Note that the Heegaard diagram constructed above might be highly inadmissible.  An
obvious periodic domain with nonnegative coefficients is given by $F$, which
represents the fibre class. However, the standard winding techniques will give us a
Heegaard diagram where $F$ (or its multiples) is the {\it only} potential periodic
domain which might prevent our Heegaard diagram from being admissible (which happens
if and only if $k=1$). In fact, we can achieve this by only changing the diagram in
the interior of $F$, so that the standard configuration of curves on
$\Sigma_{\text{min}}$ is preserved. Furthermore, we will make sure that, in the new
Heegaard diagram, the points $w$, $p^+$ and $q^+$ remain in the same connected
component. To get started, fix an arc $\delta$ in $F$, disjoint from all the $\alpha$
and $\beta$ curves and arcs in $\text{Int}(F)$, that connects the two boundary
components of $F$ and passes through $p^+$ and $q^+$. We claim that there are $g+k$
simple closed curves $\{\gamma_1,\ldots,\gamma_{g+k}\}$ in $F$ such that $\gamma_i$ do
not intersect $\delta$ and the algebraic intersection of $\gamma_i $ with $\alpha_j$
is $1$ if $i=j$ and $0$ otherwise (Note that we do not require the curves
$\gamma_1,\ldots,\gamma_{g+k}$ to be disjoint). For that, we will show that the curves
$\alpha_1,\ldots,\alpha_{g-k},\xi_1,\ldots,\xi_{2k}, \delta $ are linearly independent
in $H_1(F,\del F)$. Then the Poincar\'e-Lefschetz duality implies the existence of the
desired simple closed curves in $F$ which do not intersect $\delta$.

\begin{lemma}

The curves $\alpha_1,\ldots,\alpha_{g-k},\xi_1,\ldots,\xi_{2k}, \delta$ are linearly independent in $H_1(F,\del F)$. 

\end{lemma}

\emph{Proof.} It suffices to show that the complement of
$\alpha_1,\ldots,\alpha_{g-k},\xi_1,\ldots,\xi_{2k}, \delta$ in $F$ is connected. Take any two points $a$, $b$ in the complement. Now use the gradient flow along the northern semi-circle to obtain $\bar{a}$ and $\bar{b}$. Also let $\bar{\delta}$ be the image of $\delta$ under the flow. Connect $\bar{a}$ and $\bar{b}$ in the complement of
$\bar{\xi}_1,\ldots,\bar{\xi}_{2k}$ in $\bar{F}$ with a path that is disjoint from
$\bar{\delta}$ (This is easy because of the standard configuration of curves in
$\bar{F}$). Now flow the connecting path back to obtain a path that connects $a$ and
$b$ in the complement of $\alpha_1,\ldots,\alpha_{g-k},\xi_1,\ldots,\xi_{2k}$.
\QED

\begin{lemma}
\label{admissible}
Given a basis of the abelian group of periodic domains in the form $F, P_1,\ldots,
P_n$, after winding the $\alpha$ curves sufficiently many times along the curves $\{
\gamma_1,\ldots,\gamma_{g+k}\}$, we can arrange that any periodic domain which is given as a linear combination of $P_i$ has both positive and negative regions on the Heegaard surface. Furthermore,
for $\s \in \mathcal{S}(Y|\Sigma_{\text{min}})$, the resulting diagram is weakly admissible if $k>1$. 

\end{lemma}

\emph{Proof.} This follows by winding (\cite{OS} Section 5) successively along the curves $\{
\gamma_1,\ldots,\gamma_{g+k}\}$ in $F$, first wind along $\gamma_1$ all the $\alpha$
curves that intersect $\gamma_1$, then wind the resulting curves around $\gamma_2$,
etc. In this way the $\alpha$ curves stay disjoint (each winding is actually a
diffeomorphism of $F$ supported near $\gamma_i$, and maps disjoint curves/arcs to
disjoint curves/arcs). Furthermore, because winding along $\gamma_i$ is a
diffeomorphism of $F$ isotopic to identity, it preserves the property that $\alpha_j$
and $\gamma_k$ have algebraic intersection numbers $1$ if $j=k$, $0$ otherwise. If we
had a periodic domain with a nontrivial boundary along $\alpha_i$, then after winding
sufficiently along $\gamma_i$, the multiplicity of some region of the periodic domain
with boundary in $\alpha_i$ becomes negative. The argument for that relies on the
observation that, since the total boundary of the periodic domain has algebraic
intersection number $0$ with $\gamma_i$, and since all the other $\alpha$ curves have
algebraic intersection number $0$, while $\alpha_i$ has nonzero algebraic
intersection, the boundary of the periodic domain must also include a $\beta$ curve
which has nonzero algebraic intersection number with $\gamma_i$.  Thus after each
winding along $\gamma_i$, the domain of the periodic domain which has boundary on
$\alpha_i$ has a region where the multiplicity is decreased. Hence after sufficiently
many windings, we can ensure that any periodic domain with boundary in one of
$\alpha_1,\ldots,\alpha_{g+k}$ has at least one negative region.

Furthermore, note that a periodic domain is uniquely determined by the part of its
boundary which is spanned by $\{[\alpha_0], \ldots, [\alpha_{g+k}]\}$.  Therefore,
given a basis $F, P_1, \ldots, P_n$, after winding sufficiently many times, we can
make sure that each $P_i$ has sufficiently large multiplicities both positive and
negative in certain regions of the Heegaard diagram where all other $P_j$'s have small
multiplicities. Thus for a periodic domain to have only positive multiplicities, it
must be of the form $mF+m_1P_1+\ldots+m_nP_n$ such that $m$ is much larger than
$|m_i|$. Then $\langle c_1(\s), mF+m_1P_1+\ldots+m_nP_n \rangle = m \langle c_1(\s), F
\rangle + \sum_{i=1}^n m_i \langle c_1(\s), P_i \rangle $ must be non-zero when $k\neq
1$ since $m \langle c_1(\s), F \rangle$ dominates the sum and $\langle c_1(\s), F
\rangle = 2 - 2k $ is non-zero. Thus the diagram can be made weakly admissible when
$k>1$.  \QED

We remark that the configuration of the curves on $\bar{F}$ is left intact.  Also, the
curve $\delta$ in $F$ has not been changed. Therefore, after winding we still have the
points $p$ and $q$ lying in the same region of the Heegaard diagram.  From now on, we
will use the notation $(\Sigma, \alpha_0,\ldots,\alpha_{g+k},
\beta_0,\ldots,\beta_{g+k},z,w)$ for this diagram, which is weakly admissible if
$k>1$. We will refer to this kind of diagrams as {\it almost admissible}. In order to
make sense of Heegaard Floer homology groups for our special Heegaard diagram in the
case when the lowest genus fibre is a torus (i.e. $k=1$), we will need to work in the
perturbed setting since the periodic domain $F$ prevents the diagram from being weakly
admissible.  However, because we have an ``almost admissible'' diagram, it suffices to
perturb only in the ``direction of the fibre class''. 

\begin{lemma} 
	Given a basis of the abelian group of periodic domains in the form $F,
	P_1,\ldots, P_n$, we can find an area form $A$ on the Heegaard surface such
	that $A([F])>0$ and $A(\text{span}\{P_1,\ldots,P_n \}) = 0 $. 
\end{lemma}

\emph{Proof.} By the previous lemma, we can arrange that any periodic domain in the
linear span of $\{P_1,\ldots,P_n\}$ has both positive and negative regions on the
Heegaard surface. The rest of the proof now follows from Farkas' lemma in the theory of convex sets. See \cite{LOT} Lemma $4.17-4.18$.
\QED

Now an area form $A$ on the Heegaard surface gives a real cohomology class $[\eta] \in
H^2(Y;\f{R})$ via the bijection between periodic domains and $H_2(Y;\f{Z})$. Namely,
set $[\eta](P) = A(P)$. Choosing a representative $\eta \in [A]$ we can consider the
perturbed Heegaard Floer homology $HF^+(Y,f,\eta)$. Since $F$ is the only periodic
domain which prevents weak admissibility (only in the case $k=1$) and $\eta([F])>0 $,
we have a well-defined group $HF^+(Y,f,\eta)$ by the following lemma :

\begin{lemma}

	Given $\x,\y \in \f{T}_\alpha \cap \f{T}_\beta$,$i,j\in \f{Z}_{\geq 0}$ and
	$r,s \in \f{R}$ there are
	only finitely many homology classes $\varphi \in \pi_2(\x,\y)$,
	with $n_z(\varphi) = i-j$ and $\eta(\varphi)=r-s$ which have positive domains.
\end{lemma}
	\emph{Proof.} Let $\varphi$ and $\psi$ be in $\pi_2(\x, \y)$ ,
	then $\varphi-\psi \in \pi_2(\x,\x)$. We can write $\varphi-\psi =
	mF+m_1P_1+\ldots+ m_nP_n + n\Sigma$. Since $n_z(\varphi)=n_z(\psi)$, we have
	$n=0$. Also since $\eta(\varphi)=\eta(\psi)$ and $\eta(F)\neq 0$ while
	$\eta(P_i)=0$, we conclude that $m=0$. Finally, since $A(P_i) = 0$ , we have
	$A(\varphi) = A(\psi)$ but then there are only finitely many nonnegative domains
	which have a fixed area.
\QED

Now, as explained in the introduction $HF^+(Y,f,\eta)$ is an invariant of
$(Y,[\eta])$, in fact it only depends on $\text{Ker}(\eta) \cap H_2(Y;Z)$, hence is
independent of the value of $\eta([F])$.

The usual invariance arguments of Heegaard Floer theory, as in \cite{OS}, imply
that $HF^+(Y,f,\eta)$ is independent of the choice of $f$ within its smooth isotopy
class. Also note
that a geometric way of choosing $\eta$ is by choosing a section $\gamma$ of $f$ (a
section of $f$ always exists) and letting $[\eta]$ be the
Poincar\'e dual of $[\gamma]$. In that case, we will write $HF^+(Y,f,\gamma)$ for this
perturbed Heegaard Floer homology group. In fact, the choice of the base points $w$
and $z$ as above gives a section of $f$. Namely, note that we have arranged
so that the image of $z$ under the flow above both the northern and the southern
semi-circles lies in the same region as $w$. The union of these two gradient flow
lines can therefore be perturbed into a section of $f$, which
we will denote by $\gamma_w$. The group $HF^+(Y,f,\gamma_w)$ will be one of the main
protagonists in this paper. The differential of this group can be made more explicit
as follows: Choose a basis of the group of periodic domains in the form
$F,P_1,\ldots,P_n$ such that $F$ is the fibre of $f$ and $P_i$ are periodic domains so
that the boundary of $P_i$ does not include $\alpha_0$ or $\beta_0$ (This can be
arranged by subtracting a multiple of $F$). Then if we choose $\eta \in PD[\gamma_w]$ we
have $\eta(\text{span}(P_1,\ldots,P_m))=0$ and $\eta(F)=n_w(F)=1 $. Therefore for any
periodic domain $P$, we have $\eta(P)=n_w(P)$. Thus there exists a function $\lambda :
\f{T}_\alpha \cap \f{T}_\beta \to \f{R}$ such that for any $\varphi \in
\pi_2(\x,\y)$, we have $\eta(\varphi) - n_w(\varphi) =
\lambda(\x)-\lambda(\y)$. Hence, we can define the differential for
$HF^+(Y,f,\gamma_w)$ as follows: 

\[ \del^+([\x,i]) = \sum_{\y}
   \sum_{\varphi\in\pi_2(\x,\y), n_z(\varphi)\leq i }
   \#\widehat{\mathcal{M}}(\varphi) t^{n_w(\varphi)} [y,i-n_z(\varphi)] \]

This yields the same homology groups as the original definition where the differential
is weighted by $t^{\eta(\varphi)}$: namely, the two chain complexes are related by rescaling
each generator $[\x,i]$ to $t^{\lambda(\x)}[\x,i]$. When we
consider $HF^+(Y,f,\gamma_w)$, we will always consider the differential above.

\subsection{Splitting the Heegaard diagram}
\label{split}

As explained in the introduction, we will only consider the spin$^c$ structures on $Y$
that satisfy the adjunction equality with respect to $\Sigma_\text{min}$; the set of
isomorphism classes of such spin$^c$ structures was denoted by
$\mathcal{S}(Y|\Sigma_{\text{min}})$.  In this section we observe that for $\s\in
\mathcal{S}(Y|\Sigma_{\text{min}})$, we obtain a nice splitting of the generators of
the Heegaard Floer complex into intersections in $F$ and $\bar{F}$. Furthermore, we prove a key lemma en route to understanding the holomorphic curves contributing to the differential. 

Let us denote by $I_\lft$ the intersection of $\alpha_1 \times \ldots \times
\alpha_{g-k}$ and $\beta_1 \times \ldots \times \beta_{g-k}$ in $\sym^{g-k}(\Sigma)$,
and by $I_\rgh$ the set of 
intersection points of $\alpha_0 \times \alpha_{g-k+1} \times \ldots \times
\alpha_{g+k}$ and $\beta_0 \times \beta_{g-k+1} \times \ldots \times \beta_{g+k}$ in
$\sym^{2k+1}(\Sigma)$ such that each intersection point lies in $\bar{F}$. Thus,
each element of $I_\rgh$ consists of one point from the set of $4k$ intersection points
of $\alpha_0$ with $\eta_1,\ldots,\eta_{2k}$, another point from the set
of $4k$ intersection points of $\beta_0$ with $\xi_1,\ldots,\xi_{2k}$ and finally
$2k-1$ points from the set of $2k$ points consisting of the intersections of
$\bar{\xi}_i$ with $\bar{\eta}_i$ for $i=1,\ldots,2l.$

We have $I_\lft \otimes I_\rgh \subset \mathbb{T}_\alpha \cap \mathbb{T}_\beta$,
where $\mathbb{T}_\alpha = \alpha_0 \times \ldots \times \alpha_{g+k}$ and
$\mathbb{T}_\beta = \beta_0 \times \ldots \times \beta_{g+k}$ are the Heegaard tori in
$\sym^{g+k+1}(\Sigma)$. Denote by $C_\lft$ and $C_\rgh$ the free $\Lambda-$modules 
generated by $I_\lft$ and $I_\rgh$ respectively.

\begin{lemma} \label{splitting} An intersection point $\x\in \mathbb{T}_\alpha
	\cap \mathbb{T}_\beta $ induces a spin$^c$ structure $\s_z(\x) \in \mathcal{S}(Y|\Sigma_{\text{min}})$ if and only if $\x \in C_\lft \otimes
	C_\rgh$.
\end{lemma}

\emph{Proof.} This follows easily from the following index formula from \cite{properties} (cf. Lemma $4.11$ in $\cite{lipshitz}$):
\[\langle c_1(\s_z(\x)), F \rangle = e(F) + 2n_{\x}(F) \]
where $n_{\x}(F)$ is the number of components of the tuple $\x$ which lie in $F$. 
Since $\s_z(\x) \in \mathcal{S}(Y|\Sigma_{\text{min}})$, we have $\langle
c_1(\s_z(\x)), F \rangle = \langle
c_1(\s_z(\x)), \Sigma_{\text{min}} \rangle = 2-2k.$ Also $e(F)= -2g $, hence
the above formula gives \[ n_{\x}(F)= 1+g-k \] which is satisfied if and
only if $\x \in C_\lft \otimes C_\rgh$.
\QED

Next, we prove an important lemma about the behaviour of holomorphic disks on
the tubular regions to the left of $\alpha_0$ and $\beta_0$. This lemma lies at
the heart of most of the arguments about the behaviour of holomorphic curves
that we are going to consider subsequently. For the purpose of the next lemma,
let $a$ and $b$ be parallel pushoffs of $\alpha_0$ and $\beta_0$ to the left
into the interior of $F$. Let us label the connected components of the domains
in the cylindrical region between $a$ and $\alpha_0$ by $a_1,\ldots,a_{4k}$ and
the cylindrical region between $b$ and $\beta_0$ by $b_1,\ldots,b_{4k}$. Choose
the labeling so that $a_1$ and $b_1$ are in the same region as the arc
$\delta$, hence $n_{a_1}=n_{b_1}=n_w$. We will adapt the set-up of Lipshitz's cylindrical reformulation of Heegaard Floer homology \cite{lipshitz}. Let us also call an almost complex structure $J$ on $\Sigma \times [0,1] \times \f{R}$ admissible if it satisfies the axioms (J1-5) of \cite{lipshitz}, and the differential is obtained via a count of $J$-holomorphic curves $u: S \to \Sigma \times [0,1] \times \f{R}$ for an admissible $J$ which satisfy the axioms (M1-6) of \cite{lipshitz}. 

\begin{lemma} \label{same multiplicity} Let $\x=\x_\lft\ten\x_\rgh$ and
	$\y=\y_\lft\ten\y_\rgh$ be in $C_\lft \ten C_\rgh$ and $A\in
	\pi_2(\x,\y)$ and let $u$ be a Maslov index $1$ holomorphic curve in the
	homology class $A$. Assume moreover that the contribution of curves in the class $A$ to the differential is non-zero. Then, \[n_w(u)=
	n_{a_1}(u)=\ldots=n_{a_{4k}}(u) = n_{b_1}(u)=\ldots=n_{b_{4k}}(u)\]
	
	Furthermore, if the projection of the image of $u$ to the Heegaard
	surface $\Sigma = F \cup_{\alpha_0 \cup \beta_0} \bar{F}$ lies entirely
	in $F$, one can find almost complex structures $j_0$ and $j_1$ on
	$\Sigma$ and an admissible almost complex structure $J$ on $\Sigma
	\times [0,1] \times \f{R}$ such that $J|_{\Sigma \times \{0 \} \times
	\f{R}} = j_0$ and $J_{\Sigma \times \{1 \} \times \f{R} }= j_1$ with
	the property that $u$ restricted to the boundary does not hit $\{a, b \} \times \{0,1\} \times \f{R}$. (In
	other words, $u$ converges to Reeb orbits around $a$ and $b$ upon neck
	stretching).  \end{lemma}

\emph{Proof.} The first part of the proof will be obtained by ``stretching the neck'' along the curves $a$ and $b$. Suppose that there is an $i$ (mod $4k$) such that $n_{a_i}(u)\neq n_{a_{i+1}}(u)$ (one
can argue in the same way for $b_i$'s). Thus the source $S$ of $u$ has a piece of boundary which maps to the $\beta$ arc that separates $a_i$ and
$a_{i+1}$. Let $\beta_j$ be the curve containing that arc. The crucial
observation that we will make use of is the fact that the disk $u$ has no corners in $\beta_j \cap F$, since $\x$ and
$\y$ have no components in $\beta_j \cap F$. 

We now degenerate $\Sigma$ along the curves $a$ and $b$. Specifically, this means that
one takes small cylindrical neighborhoods of the curves $a$ and $b$, and changes the
complex structure in that neighborhood so that the modulus of the cylindrical
neighborhoods gets larger and larger. Topologically this degeneration can be
understood as follows: After degenerating along $a$ and $b$, $\Sigma$
degenerates into $\Sigma_{\text{max}}$ and $\Sigma_{\text{min}}$ and the homology
class $A$ splits into $A_\lft$ and $A_\rgh$ corresponding to the induced domains on
$\Sigma_{\text{max}}$ and $\Sigma_{\text{min}}$ from the domain of $A$ on $\Sigma$.
(The definition of homology classes $\pi_2(\x,\y)$ in this degenerated setting is
given in Definition 4.14 of \cite{LOT}. It is the homology classes of maps to $
\Sigma_{\text{max}} \times [0,1] \times \f{R}$ (and to $\Sigma_{\text{min}} \times
[0,1] \times \f{R})$ which have strip-like ends converging to $\x$ and $\y$, and to
Reeb chords at points of degeneration). 

Next we analyze the holomorphic degeneration of $u$. Suppose that the moduli
space of holomorphic curves representing $A$ is non-empty for all large values
of the stretching parameter. Then we conclude by Gromov compactness that there
is a subsequence converging to a pair of holomorphic combs of height $1$ (in
the sense of \cite{LOT} Section 5.4; see Proposition 5.23 for the proof of
Gromov compactness in this setting) $u_0$ representing $A_\lft$ and $u_1$
representing $A_\rgh$. (The limiting curves have height $1$ because otherwise
one of the stages would have index $\leq 0$, contradicting transversality --
see Proposition 5.6 of \cite{LOT}). By assumption, the degeneration of $u$
involves breaking along a Reeb chord $\rho$ contained in $a$ with one of the
ends of $\rho$ on $a \cap \beta_j $. Hence some component $S_0$ of the domain of $u_0$
has a boundary component $\Gamma$, consisting of arc components separated by
boundary marked points, such that one of the arcs is mapping to $\beta_j$ and,
at one end of that arc, $u_0$ has a strip-like end converging to the Reeb chord
$\rho$. Now, since there are no corner points on any of the $\beta$-arcs in
$\Sigma_{\text{max}}$, the marked points on $\Gamma$ are all labeled by Reeb
chords on $a$ (corresponding to arcs connecting intersection points of $\beta$
curves with $a$), and any two consecutive punctures on $\Gamma$ are connected
by an arc which is mapped to part of a $\beta$ arc which lies on the left half
of the Heegaard diagram. Thus, in particular there are no arcs in $\Gamma$
which map to $\alpha$ curves. Now, we can extend $u_0$ at the punctures on
$\Gamma$ by sending the marked points to the point of $\Sigma_{\text{max}}$ to
which $a$ collapses upon neck-stretching (This is possible since, after
collapsing $a$, $u_0|_{S_0}$ viewed as a map to $\Sigma_{\text{max}}$ admits a
continuous extension at these points. Note that the projection to $[0,1] \times
\f{R}$ also extends continuously at the punctures by the definition of
holomorphic combs; see the proof of Proposition 5.23 of \cite{LOT} for more
details regarding this).  Therefore, the image of the boundary component
$\Gamma$ under the projection to $[0,1]\times \f{R}$ remains bounded and is
entirely contained in ${0}\times \f{R}$. Moreover, since the projection is
holomorphic, the projection of $\Gamma$ to $0\times \f{R}$ is a non-increasing
function. Hence we conclude that $\Gamma$ maps to a constant.  Now, the
maximum principle implies that the entire component $S_0$ has to be mapped to a
constant value in $0 \times \f{R}$. Therefore, $S_0$ has all of its boundary
components mapped to $\beta$ curves. Now, $u_0$ maps all of its boundary to
$\beta$ curves in $\Sigma_{\text{max}}$ which remain linearly independent in
homology even after degeneration. However, $u_0(S_0)$ gives a homological
relation between those curves. The only way this could be is if the this
relation is trivial, that is, the boundary of $u_0(S_0)$ traces each $\beta$
curve algebraically zero times, but that contradicts the assumption that
$n_{a_i}({u_0}_{|S_0})\neq n_{a_{i+1}}({u_0}_{|S_0})$ and thus proves the first
part of the lemma. 

To prove the second part, suppose that $u$ contributes to the differential
between $\x=\x_\lft \ten \x_\rgh$ and $\y = \y_\lft \ten \y_\rgh $ and the
image of $u$ lies entirely in $F$ (the left side of the Heegaard diagram). This
implies that $\x_\rgh = \y_\rgh$ by using the first part of the lemma which we
have already established. Let us describe how we choose the complex structure
$j_0$ ($j_1$ is constructed in a completely analogous way). Let $\beta_j$ be
the curve such that $\beta_j \cap \alpha_0$ is an intersection point that
appears in $\x_\rgh = \y_\rgh$. Note that $u$ cannot have any boundary
component that maps to any other $\beta$ curves that intersect with $\alpha_0$.
Recall also that we have the closed $\beta$ curves $\beta_1, \ldots,
\beta_{g-k}$ in $F$.  After stretching the neck around $a$ and $b$
sufficiently, suppose that $u$ restricted to the boundary still intersects $a$.
As before, in the limit $u$ degenerates and we restrict our attention to the
component $u_0$ which has boundary component that maps to $\beta_j$. We
identify the left side of the degenerated Heegaard surface with
$\Sigma_\text{max}$. From now on, we also think of $\beta_j$ as a closed curve
since after the degeneration along $a$, the two end points of $\beta_j$ come
together. By exactly the same argument as in the first part, we conclude that
$u_0$ maps all of its boundary components to $\beta_1 \cup \ldots \cup
\beta_{g-k} \cup \beta_j$ and its projection to $[0,1] \times \f{R}$ is
constant and lies on $\{0\} \times \f{R}$. To arrive at a contradiction, we
would like to ensure that no such $u_0$ exists by choosing the almost complex
structure $J$ on $\Sigma \times [0,1] \times \f{R}$ appropriately such that the
resriction of the induced complex structure complex structure on
$\Sigma_{\text{max}} \times [0,1] \times\f{R}$ to ${\Sigma_{\text{max}} \times
\{0 \} \times \f{R}}$ does not allow such a curve.  To that end, we adopt the
idea used in Lemma 8.2 of \cite{lipshitz} (which in turn is adapted from the
idea in Proposition 3.16 of \cite{OS}). Namely, since all of $\beta_1, \ldots,
\beta_{g-k}, \beta_j$ are linearly independent, we can find disjoint curves on
$\Sigma_{\text{max}}$ not intersecting these $\beta$ curves such that their
complement in $\Sigma_{\text{max}}$ is a disjoint union of $(g-k+1)$ punctured
surfaces with each surface having genus at least $1$ and such that each $\beta$
curve is contained in one and only one of these surfaces. We further degenerate the
complex structure by stretching along these curves. Note that, crucially, these
curves are also disjoint from the original curves $a$ and $b$ along which we
degenerate. Therefore we can do the degeneration simultaneously.  Let us denote
by $j_n$ a sequence of complex structures on $\Sigma_{\text{max}}$ where we
stretch along the specified curves as $n$ tends to infinity. Following the
argument in Lemma 8.2 of \cite{lipshitz}, suppose we have a degeneration $J_n$
of admissible almost complex structures on $\Sigma \times [0,1] \times \f{R}$,
corresponding to stretching along $a$ and $b$ and when restricted to a
neighborhood of $\Sigma \times \{0\} \times \f{R}$ , it has a further
degeneration of the form $J_n= j_n \times j_{[0,1]\times \f{R}}$ corresponding
to stretching along the other curves that we dicussed (which we said are
disjoint from $a$ and $b$). Suppose now that there is a sequence of $J_n$
holomorphic curves  that converge to a curve whose projection to $[0,1] \times
\f{R}$ is constant to a point $p \in \{0 \} \times \f{R}$. By composing $u_n$
with projection to $[0,1] \times \f{R}$, and rescaling near $p$, we get a
sequence of maps $u_n: (S_n, \partial S_n)  \to (\f{H}, \f{R})$, where $S_n$ is
obtained from $\Sigma_{\text{max}}$ by cutting along the $\beta$ curves. On the
other hand, by hypothesis the projection of $u_n$ to $\Sigma$ is supported in
$F$, which in the limit degenerates to $g-k+1$ disjoint surfaces of genus at
least 1 such that each $\beta$ curve is in a separate surface. Therefore, by
passing to a subsequence and restricting to the complement of the $\beta$
curves, in the limit we obtain a $(g-k+1)$-fold covering map from a disjoint
union of $g-k+1$ punctured surfaces of genus at least 1 to $(\f{H},\f{R})$,
which cannot exist. Hence, for a sufficiently large $n$, if we set $j_0 =
j_n$, we can conclude that the map $u_0$ can not exist. Hence, by stretching as
described, we can ensure that $u$ does not map any of its boundary to $a$ for
some $J_n$ where $n$ is sufficiently large. One can similarly arrange by a
stretching supported near $\Sigma \times \{1 \} \times \f{R}$  so that $u$
restricted to its boundary does not intersect $b$. \QED

\subsection{Calculations for fibred $3$-manifolds and $C_\rgh$}

Before delving into a general study of Heegaard Floer homology for broken
fibrations, here we will calculate $HF^+(Y,\eta)$ in the case of fibred
$3$--manifolds. Some of these calculations were done independently by Wu in
$\cite{Wu}$, where perturbed Heegaard Floer homology for $\Sigma_g \times S^1$
is calculated for all spin$^c$ structures. We take the liberty to reconstruct
some of the arguments presented there in this section since these calculations
will play a role for the calculations we do for general fibred $3$--manifolds.
Even though we will do calculations in general for any fibred $3$--manifold, we
will restrict ourselves to spin$^c$ structures in $\mathcal{S}(Y|F)$, which
will simplify the calculations. Our conclusion is that $\oplus_{\s\in
\mathcal{S}(Y|F)}HF^+(Y,\eta)$ has rank $1$. See also \cite{AiPeters} for a
different approach in the case of torus bundles.  

For fibred $3$--manifolds, we have $g=k$, thus the Heegaard diagram has the
curves $\alpha_0$,$\beta_0$, and the rest of the diagram is constructed from
the standard configuration of curves $\bar{\xi}_1,\ldots,\bar{\xi}_{2k}$,
$\bar{\eta}_1,\ldots,\bar{\eta}_{2k}$ as in Figure \ref{curves}. Also we will
see below that, for the spin$^c$ structures in $\mathcal{S}(Y|F)$, the
generators of our chain complex are given by the intersection points in
$C_\rgh$. 

We first discuss the case of torus bundles. It will then be clear that the
general case is just a matter of notational complication. Also note that, in
the case of torus bundles, we have to use a perturbation $\eta$ with $\eta([F])
> 0$ as explained in the previous section since our diagram is not weakly
admissible. For higher genus fibrations, the diagram is weakly admissible hence
our calculation also determines the unperturbed Heegaard Floer homology
$HF^+(Y)$. When doing explicit calculations we will always consider the case of
$HF^+(Y,f,\gamma_w)$ but clearly all arguments go through for any perturbation
with $\eta$ satisfying $\eta([F])>0$, or for the unperturbed case whenever the
diagram is weakly admissible. 

\begin{figure}[!h] \centering \includegraphics[scale=0.8]{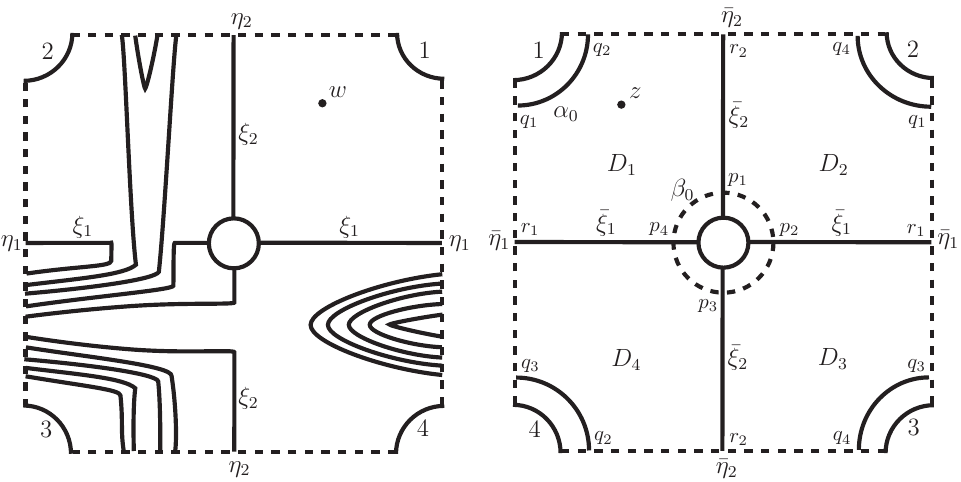}
	\caption{Torus bundles} \label{torus bundles} \end{figure}

Figure \ref{torus bundles} shows the Heegaard diagram for $T^3$. Both the left
and the right figure are twice punctured tori, and are identified along the two
boundaries (the one in the middle and the one formed by the four corners) where
the gluing of the left and right figures is made precise by the labels at the
four corners. On the right side the standard set of arcs
$\bar{\xi}_1,\bar{\xi}_2, \bar{\eta}_1, \bar{\eta}_2$ are depicted; the left
side is constructed by taking the images of these arcs under the horizontal
flow (which is the identity map for $T^3$), and winding $\xi_1$ and $\xi_2$
along transverse circles so that the diagram becomes almost admissible (Note
that the winding process avoids the region where $w$ is placed, as required:
first $\xi_2$ is wound once along a horizontal circle, then $\xi_1$ is wound
twice along a vertical circle). For general torus bundles, the same
construction will give a Heegaard diagram, where $\xi_1$ and $\xi_2$ are
replaced by their images under the monodromy of the torus bundle.  The
important observation here is that the right side of the diagram is always
standard. We will show that all the calculations that we need can be done on
the right side of the diagram for the spin$^c$ structures we have in mind. The
calculation for $T^3$ is essentially the same as in \cite{Wu}. However, we will
see that Lemma \ref{same multiplicity} plays a crucial role in the calculation
for general torus bundles. We first do the calculation for $T^3$.

\begin{proposition} \label{3torus}  $HF^+(T^3, f, \gamma_w, \s_0) = \Lambda $
	where $\s_0 \in \mathcal{S}(T^3 | T^2)$ is the unique torsion spin$^c$
	structure on $T^3$.  \end{proposition}

\emph{Proof.} As in Lemma \ref{splitting}, $\s_\x(z) \in \mathcal{S}(T^3|T^2)$
if and only if $\x \in C_\rgh$, hence $\x$ can be one of the following tuples
of intersections depicted in Figure \ref{torus bundles}: \[ \x_1 =  p_1q_2r_1 \
\ \x_2 = p_2q_1r_2 \ \x_3 = p_3q_4r_1 \ \ \x_4 = p_4q_3r_2 \] \[ \y_1 =
p_4q_1r_2 \ \ \y_2 = p_1q_4r_1 \ \ \y_3 = p_2q_3r_2 \ \ \y_4 =  p_3q_2r_1 \]

Next, we apply the adjunction inequality for the other $T^2$ components, this
implies that the Heegaard Floer groups vanish except for the unique torsion
spin$^c$ structure, $\s_0$ which has $c_1(\s_0)=0$. The two other torus
components are realized by periodic domains in Figure \ref{torus bundles} , one
of them is the domain $P_1$ including $D_2\cup D_3$ and bounded by $\alpha_2$
and $\beta_1$, the other one is the domain $P_2$ including $D_3 \cup D_4$ and
bounded by $\alpha_1$ and $\beta_2$. Then the formula $\langle c_1(\s_{z}(\x)),
P_i \rangle = e(P_i) + 2n_\x(P_i)$, implies that the only intersection points
for which $\s_z(\x)=\s_0$ are $\x_1$ and $\y_1$. Furthermore, note that $D_1$
is a hexagonal region connecting $\x_1$ to $\y_1$, hence it is represented by a
holomorphic disk $\varphi_1 \in \pi_2(\x_1,\y_1)$, and the algebraic number of
holomorphic disks in the corresponding moduli space of disks in the homology
class of $\varphi_1$ is given by $\#\widehat{\mathcal{M}}(\varphi_1)=\pm 1$
(See appendix in \cite{Rasmussen}).

Now, given any other Maslov index $1$ homology class $A \in \pi_2(\x_1,\y_1)$,
we have $A = D_1 + mF+  m_1 P_1 + m_2 P_2$. In particular, note that
$n_z(A)=1$.  Furthermore, if we restrict to those with $n_w = 0$ (that is
$m=0$), since $m_1 P_1 + m_2 P_2$ has both positive and negative domains by
almost admissibility, unless $m_1=m_2=0$ there is no holomorphic representative
of $A$.  

We conclude that $\del^+[\x_1,i] = f(t) [\y_1,i-1]$, where $f(t)=\pm 1 +
\OO{t}$ is invertible in the Novikov ring. This implies that $[\y_1,i]$ is in
the image of $\del^+$. Thus in particular we have $\del^+[\y_1,i]=0$ for all
$i$. Finally, there is no Maslov index $1$ disk with $n_w=0$ which connects
$\x_1$ to itself or $\y_1$ to itself. Thus we conclude that in
$CF^+(T^3,f,\gamma_w, \s_0)$: \[\del^+[\x_1,0]=0\ \ \ \del^+[\y_1,i]=0 \] \[
\del^+[\x_1,i] = (\pm1 + \OO{t})[\y_1,i-1] \ \ \ \text{for} \ i>0  \] Hence the
homology is generated by $[\x_1,0]$, in other words $HF^+(T^3, f, \gamma_w,
\s_0) = \Lambda$ as required. \QED

From now on, we will simply write $\x_1$ for $[\x_1,0]$. The next theorem
generalizes this calculation to any torus bundle.

\begin{theorem} \label{toruscalc} Let $(Y,f)$ be a torus bundle and let $\s$ be
	in $\mathcal{S}(Y | T^2)$. Then, $HF^+(Y, f, \gamma_w,\s) = \Lambda $
	if $\s = \s_0$ where $\s_0$ is the spin$^c$ structure corresponding to
	vertical tangent bundle and $HF^+(Y,f,\gamma_w,\s)= 0$ otherwise.
\end{theorem}

\emph{Proof.} The main difficulty for the general torus bundle case that makes
the calculation different from the calculation for $T^3$ is that we cannot a
priori eliminate the generators $\x_2,\x_3,\x_4$ and $\y_2,\y_3,\y_4$. In fact,
if the first Betti number of the torus bundle is equal to $1$, these generators
are in the same spin$^c$ class as $\x_1$ and $\y_1$.

Now, the domains $D_i$ are homology classes in $\pi_2(\x_i,\y_i)$, which have
holomorphic representatives $\varphi_i$ with $\# \mathcal{M}(\varphi_i) = \pm
1$.  Since any non-trivial periodic domain has to pass through some region to
the left of $\alpha_0$ or $\beta_0$, any other homology class in
$\pi_2(\x_i,\y_i)$ which contributes to the differential has to have $n_w \neq
0$ by Lemma \ref{same multiplicity}. For the same reason, any homology class in
$\pi_2(\x_i, \y_j)$ for some $i \neq j$ which contributes to the differential
has to have $n_w \neq 0$ since there is no homology class in $\pi_2(\x_i,\y_j)$
that lies in the right side of the diagram (this can be verified either by
inspection, or referring to the case of $T^3$, where $\x_i$ and $\y_j$
represent different spin$^c$ classes for $i\neq j$). Moreover, the classes in
$\pi_2(\x_i,\x_j)$ all have even Maslov index, hence do not contribute to the
differential. Therefore, we have \[ \del^+[\x_1, i] = [\y_1,i-1] \ \ \
(\text{mod}\ t) \ \ \text{for} \ \ i>0 \] \[ \del^+[\x_1,0]= 0 \ \ \
(\text{mod}\ t) \ \ \ \del^+[\x_2, i] = [\y_2,i] \ \ \  (\text{mod}\ t) \]  \[
\del^+[\x_3, i] = [\y_3,i] \ \ \  (\text{mod}\ t) \ \ \ \del^+[\x_4, i] =
[\y_4,i] \ \ \  (\text{mod}\ t)\]

where the higher order terms do not involve the $\x_j$'s. As before, because we
are working over a Novikov ring of power series, we conclude that
$[\y_1,i],[\y_2,i],[\y_3,i]$ and $[\y_4,i]$ are all in the image of $\del^+$.
Furthermore, the only possible generator which might be in the kernel of
$\del^+$ is $[\x_1,0]$. Finally lemma \ref{no curve} below shows that there is
no holomorphic disk starting at $\x_1$ with $n_z=0$ and $n_w\neq 0$. Hence we
have $\del^+[\x_1,0]=0$ and the homology group
$\oplus_{\s\in\mathcal{S}(Y|T^2)} HF^+(Y,f,\gamma_w,\s)$ is generated by
$[\x_1,0]$. Furthermore, $\s_z(\x_1)=\s_0$ so the theorem is proved.  \QED

Note also that the adjunction inequality implies that $HF^+(Y, f, \gamma_w, s)$
vanishes for $\s \not \in \mathcal{S}(Y|T^2)$. Therefore the above calculation
is in fact a complete calculation of perturbed Heegaard Floer homology for
torus bundles.

The following lemma which we used in the above calculation holds in general
(not only in the fibred case). Let $Y$ be any $3$--manifold with $b_1>0$, and
$f:Y \to S^1$ a broken fibration with connected fibres.  Construct the almost
admissible Heegaard diagram for $f$ as before and let $\x_1 \in C_\rgh$ be
given by the union of the intersection points in $\alpha_0 \cap \bar{\eta}_2 $,
$\bar{\xi}_2 \cap \beta_0$ , and $\bar{\xi}_i \cap \bar{\eta}_i $ for $i \neq
2$, where the intersection point in $\alpha_0 \cap \beta_2$ and $\alpha_2 \cap
\beta_0$ are chosen so that the region containing $z$ includes them as corners.
(In the case of the torus bundle this is the generator $[\x_1, 0]$). Note that
the generators of $C_\rgh$ can always be described from the standard diagram
since the right hand side of our Heegaard diagrams is always the same. 

\begin{lemma} \label{no curve} Let $\varphi \in \pi_2(\x_\lft \ten \x_1,
	\y_\lft \ten \y_\rgh)$ be a holomorphic disk in a class that
	contributes non-trivially to the differential for given $\x_\lft,
	\y_\lft, \y_\rgh$. If $n_z(\varphi)= 0$, then $\y_\rgh = \x_1$ and the
	domain of $\varphi$ is contained on the left side of the Heegaard
	diagram (i.e. it is contained in $F$). 

\end{lemma}

\emph{Proof.} Consider the component of $\x_1$ which is an intersection point
on $\beta_0$, say $p_1$. Now, among the four regions which have $p_1$ as one of
their corners, one includes $z$, namely $D_1$, and two of them lie in the left
half of the diagram, hence by lemma \ref{same multiplicity}, they must have the
same multiplicity. Denote these regions by $L_1$ and $L_2$, so that $L_1$ and
$D_1$ share an edge on $\beta_0$. If the component of $\varphi$ which is
asymptotic to $p_1$ is constant, then $p_1$ is also part of $\y_\rgh$.
Otherwise, since $\varphi$ has a corner which leaves $p_1$ and
$n_z(\varphi)=0$, we must have a non-zero multiplicity at $L_2$, but since
$L_1$ and $L_2$ must have the same multiplicity, this implies that $p_1$ has to
be a member in $\y_\rgh$. The same conclusion applies for the point of $\x_1$
which lies on  $\alpha_0$. But then there is a unique way to complete these two
intersection points to a generator in $C_\rgh$, hence we conclude that $\y_\rgh
= \x_1$. Thus $\varphi$ is in $\pi_2(\x_\lft\ten\x_1,\y_\lft\ten\x_1)$. 

Furthermore, since $\varphi$ fixes $\x_1$, the intersection of the domain of
$\varphi$ with $\bar{F}$ must coincide with the intersection of some periodic
domain for $S^1\times \Sigma_k$ with $\bar{F}$ (since any domain that has no
corners on the right side, can be completed to a periodic domain on the
Heegaard diagram of $S^1 \times \Sigma_k$ by reflecting). However, it is easy
to identify all the periodic domains of $S^1\times \Sigma_k$ and observe that
no non-trivial combination of periodic domains for $S^1\times \Sigma_k$ (if we
leave out $F$ and its multiples), can have the same multiplicity in the regions
immediately to the left of $\alpha_0$ and $\beta_0$.  However, by Lemma
\ref{same multiplicity} this property has to hold. This proves the lemma. \QED 

\begin{theorem} \label{fibred} Let $(Y,f)$ be a fibre bundle with fibre a genus
	$g$ surface and let $\s$ be in $\mathcal{S}(Y | \Sigma_g)$. Then,
	$HF^+(Y, f, \gamma_w,\s) = \Lambda $ if $\s = \s_0$ where $\s_0$ is the
	spin$^c$ structure corresponding to vertical tangent bundle and
	$HF^+(Y,f,\gamma_w,\s)= 0$ otherwise.  \end{theorem} \emph{Proof.} The
	proof is essentially the same as the proof of the corresponding theorem
	for the torus bundles. The only difference is the number of generators
	which are cancelled out: there are now $8g$ generators $\x_1,\ldots
	\x_{4g}$, $\y_1,\ldots,\y_{4g}$, and the $4g$ hexagonal regions of
	$\bar{F}$ (see Figure \ref{curves}) give $\del^+[\x_1, i] = [\y_1,i-1]
	\ \  (\text{mod}\ t)$ and $\del^+[\x_j,i] = [\y_j,i] \ \ (\text{mod}\
	t)$ for $j\geq 2$. Arguing as before, the only generator left is again
	$\x_1$ which gives $\s_z(\x_1) =\s_0$. \QED

Note that this gives a new way of obtaining the results of the original
calculation of Ozsv\'ath and Szab\'o in \cite{OS} for fibred $3$--manifolds. 

\begin{corollary} Let $(Y,f)$ be a fibre bundle with fibre a genus $g>1$
	surface and let $\s$ be in $\mathcal{S}(Y | \Sigma_g)$. Then, $HF^+(Y,
	\s) = \f{Z} $ if $\s = \s_0$ where $\s_0$ is the spin$^c$ structure
	corresponding to vertical tangent bundle and $HF^+(Y,\s)= 0$ otherwise.

\end{corollary} \emph{Proof.} Since the diagram is weakly admissible, we can
let $t=1$ and the result follows from the previous theorem . \QED

In general, let $\del^+_\rgh$ be the contribution to the Heegaard Floer
differential from the holomorphic disks whose domain lies in $\bar{F}$ (i.e.
the disks which lie on the right half of our almost admissible Heegaard
diagrams), also we write $CF^+_\rgh = C_\rgh \ten \Lambda[\f{Z}_{\geq 0}]$ for
the chain complex associated with the right side of the diagram for the purpose
of constructing $HF^+$ theory, that is, the chain complex freely generated over
$\Lambda$ by $[\mathbf{x},i]$ where $\mathbf{x} \in C_\rgh$ and $i \in
\f{Z}_{\geq 0}$.  \begin{corollary} \label{Crgh} $(CF^+_\rgh,\del^+_\rgh)$ is a
	chain complex with rank $1$ homology generated by $\x_1$.
\end{corollary} \emph{Proof.} This is only a reformulation of the above
results. \QED

\section{The isomorphism}
\label{section3}
In this section, we prove the main theorem of this paper. Namely, we prove that the
perturbed Heegaard Floer homology group $HF^+(Y,f,\gamma_w)$ is isomorphic to the
Floer homology of the chain complex $(C_\lft, \del_\lft; \Lambda)$. Before stating
our theorem let us digress to give a rigorous definition of the latter chain complex. 

\subsection{A variant of Heegaard Floer homology for broken fibrations over the
circle} \label{variant} Let $Y$ be a 3--manifold with $b^1>0$, and let $f:Y\to S^1$ be
a broken fibration with connected fibres, and satisfying the conditions at the
beginning of Section \ref{section2.1}. As before consider the highest genus fibre
$\Sigma_g$ and let $\alpha_1,\ldots, \alpha_{g-k}$ and $\beta_1,\ldots,\beta_{g-k}$ be
tuples of $g-k$ disjoint linearly independent simple closed curves on $\Sigma_g$
obtained from the attaching circles corresponding to the critical values of $f$, and let $w$ be a base point that is in the complement of $\alpha$ and $\beta$ curves. As in Lemma
\ref{admissible}, we can arrange weak admissibility for $k>1$ by winding if necessary. For $k=1$, we need to keep track of the intersection with the point $w$ and have to work over $\Lambda$. We define the Floer homology of such a configuration in a manner similar to
the usual Heegaard Floer theory by defining the chain complex to be the
$\Lambda-$module freely generated by intersection points of $\f{T}^{g-k}_\alpha =
\alpha_1 \times \ldots \alpha_{g-k}$ and $\f{T}^{g-k}_\beta =\beta_1 \times \ldots
\beta_{g-k}$ in $\sym^{g-k}(\Sigma_{g})$, equipped with a differential given as follows: \[ \del \x = \sum_{\varphi\in\pi_2(\x,\y),
\mu(\varphi)=1} \#\widehat{\mathcal{M}}(\varphi) t^{n_w(\varphi)} \y \] 

For reasons that will be clarified in Section \ref{section4}, we will denote
the homology group that we expect to get from this construction
$QFH'(Y,f;\Lambda)$. This stands for {\it quilted Floer homology} of the broken
fibration $(Y,f)$ with coefficients in $\Lambda$. There are at least two
obvious issues that we need to address in order to make sure that
$QFH'(Y,f;\Lambda)$ is well-defined. The first issue is the compactness of the
moduli space $\mathcal{M}(\varphi)$. The second issue is proving that
$\del^2=0$. The setup here is more delicate than the usual setup of Heegaard
Floer homology due to the fact that $\sym^{g-k}(\Sigma_g)$ is not a
(positively) monotone symplectic manifold when $k>0$ (it has $\langle c_1,
[\Sigma_g]\rangle=2-2k$).  Therefore, one expects the existence of
configurations with negative Chern number bubbles. However, we will adopt the
cylindrical setting of Lipshitz (\cite{lipshitz}), whereby one considers
pseudo-holomorphic curves in $\Sigma_g\times [0,1]\times \f{R}$ instead of
disks in $\sym^{g-k}(\Sigma_g)$, and choose our almost complex structures from
a sufficiently general class.  Namely, one chooses a translation-invariant
almost-complex structure $J$ on $\Sigma_{g} \times[0,1]\times\f{R}$ such that
$J$ preserves a $2$-plane distribution $\xi$ on $\Sigma_g \times [0,1]$ which
is tangent to $\Sigma_g$ near $(\boldsymbol{\alpha}\cup \boldsymbol{\beta})
\times [0,1]$ and near $\Sigma_g \times \del [0,1]$ (see \cite{lipshitz}, axiom
J5' ). Now we can show that transversality can be achieved for holomorphic
curves in the homology class of the fibre of the projection $\Sigma_g \times
[0,1]\times \f{R} \to [0,1]\times \f{R}$. However the expected dimension of
these curves is negative, hence bubbling  at interior points can be ruled out a
priori (see \cite{lipshitz} Lemma 8.2). Furthermore, since we assumed that all
the fibres are connected, the $(g-k)$-tuples of curves are linearly independent
in homology; this implies that any boundary bubble lifts to a spherical class
in $\pi_2(\sym^{g-k}(\Sigma_g))$. By choosing almost complex structures such
that $J_{|\Sigma\times \{0\} \times\f{R}}$ and $J_{|\Sigma\times \{1\} \times
\f{R}}$ are of special type as described in Lemma 8.2 of \cite{lipshitz} (cf.
Lemma \ref{same multiplicity}), we can also avoid boundary bubbles. Hence, the
compactness of $\mathcal{M}(\varphi)$ is ensured. We will once and for all
restrict our choice of almost complex structures to this class of almost
complex structures.

The drawback of this approach is that it does not correspond in a straightforward way
to the original setting in $(\sym^{g-k}\Sigma_{g},\f{T}^{g-k}_\alpha,
\f{T}^{g-k}_\beta)$ since such general almost complex structures prevent the fibres of
the projection to $[0,1] \times\f{R}$ from being complex. In this case, in order to be
able work in $\sym^{g-k}(\Sigma_g)$ one needs to establish
a proper combinatorial rule for handling bubbled configurations (for example by
applying the general machinery of virtual fundamental cycles \cite{liuTian}). It is
reasonable to expect that one would then get the same differential as above, but the
argument would be technically very involved. However, there is an
exception to this, namely when we are in the strongly negative case, that is when
$g<2k$. We show in Section \ref{section4} that in this case we can indeed use
integrable complex structures of the form $\sym^{g-k}(j_s)$ for a path $j_s$ of
complex structures on $\Sigma$ and still avoid bubbling since such complex structures are sufficient to achieve transversality and the assumption $g<2k$ ensures that the expected dimension of bubbles is negative. 

The proof of $\del^2=0$ for $QFH'(Y,f;\Lambda)$ will be part of the proof of the
isomorphism that we will construct between $QFH'(Y,f;\Lambda)$ and
$HF^+(Y,f,\gamma_w)$.  Namely, this follows from an identification between the Maslov
index $1$ moduli spaces in both theories. Furthermore, we will also see in this
section that $QFH'(Y,f;\Lambda)$ is an invariant of $(Y,[f])$, that is it only depends
on $Y$ and the homotopy class of $f$ through the homology class $[\Sigma_\text{min}] \in H_2(Y)$.

As usual in Floer homology theories, the groups $QFH'(Y,f;\Lambda)$ are graded by equivalence
classes of spin$^c$ structures. Given an intersection point in $\x\in
\f{T}^{g-k}_\alpha \cap \f{T}^{g-k}_\beta$ one gets a spin$^c$ structure $\s(\x) \in
\mathcal{S}(Y|\Sigma_{\text{min}})$, as
in Heegaard Floer theory, except we do not need to consider any additional base point
since the intersection point $\x$ gives a matching of index 1 and 2 critical points of
$f$, which in turn determines a spin$^c$ structure by taking the gradient vector field
of $f$ outside of tubular neighborhoods of these matching flow lines and extending it
in a non-vanishing way to the tubular neighborhoods.  We remark that in our setup of
Heegaard diagram for $(Y,f)$, we have the equality $\s(\x_\lft) = \s_z(\x_\lft \ten
\x_1)$ (where $\x_1$ is as in Lemma \ref{no curve}). 

{\bf Remark:} Note that if we restrict to the case where we only count $n_w=0$ curves we obtain Juh\'asz's sutured Floer homology groups associated
with the diagram $(F,\alpha_1,\ldots,\alpha_{g-k},\beta_1,\ldots,\beta_{g-k})$ (see
\cite{juhasz}). We will return to this below.

\subsection{Isomorphism between $QFH'(Y,f;\Lambda)$ and $HF^+(Y,f,\gamma_w)$}

We now proceed to prove an isomorphism between $QFH'(Y,f;\Lambda)$ and
$HF^+(Y,f,\gamma_w)$. As a first step, we make use of the calculations of the previous
section. Let $CF^+_\lft=C_\lft \ten \Lambda[\f{Z}_{\geq 0}] $ and $CF^+_\rgh = C_\rgh
\ten \Lambda[\f{Z}_{\geq 0}]$, using the splitting of generators of
$HF^+(Y,f,\gamma_w)$ as discussed in Section \ref{split}, so that we have
$CF^+(Y,f,\gamma_w)= CF^+_\lft \ten CF^+_\rgh$.  We denote by $\del_{F}$ and
$\del_{\bar{F}}=\mathds{1}\ten \del_\rgh$ the contributions to the Heegaard Floer
differential from holomorphic curves whose domains lie in $F$ and $\bar{F}$
respectively. Furthermore, we denote by $\del_\lft \ten \mathds{1}$, the contribution
of those holomorphic curves whose domain lies in $F$ and which act by identity on
$C_\rgh$ with respect to the splitting $C_\lft\ten C_\rgh$. Since
the boundary of $F$ includes points of intersections occurring in $C_\rgh$, this is a
priori more restrictive than $\del_{F}$. However, Lemma \ref{no curve} implies that $\del_\lft\ten \mathds{1}$ is a differential on $C_\lft\ten \x_1$ (and the argument given in the proof of Lemma \ref{no curve} more generally shows that $\del_F = \del_\lft \ten \mathds{1}$). The next proposition says that the homology of this differential is isomorphic to
$HF^+(Y,f,\gamma_w)$. 

\begin{proposition} \label{prop1}

	$ HF^+(Y, f,\gamma_w, \s ) \simeq H(C_\lft \ten \x_1 , \del_\lft\ten
	\mathds{1}, \gamma_w, \s)$ for $\s\in \mathcal{S}(Y|\Sigma_\text{min})$.

\end{proposition}

\emph{Proof.} Both homology groups are filtered by $n_w$. Therefore, there are
spectral sequences converging to both sides induced by the $n_w$ filtration.
Furthermore, we claim that there is a chain map: \[ F:C_\lft \ten \x_1 \to CF^+_\lft
\ten CF^+_\rgh \] given by \[ F(\x_\lft \ten \x_1) =  [\x_\lft \ten \x_1, 0] \] which
induces an isomorphism of $E^1$--pages of the spectral sequences associated with both
chain complexes. The fact that $F$ is a chain map, is a consequence of Lemma
\ref{no curve}. More precisely, Lemma \ref{no curve} gives that if a holomorphic map
contributing to the differential originates at  $[\x_\lft\ten\x_1,0]$ then it has to
converge to a generator of the form $[\y_\lft\ten\x_1,0]$, and the domain of the map
has to lie on the left half of the Heegaard diagram; these are
exactly the contributions to the differential captured by $\del_\lft \ten \mathds{1}$.

Furthermore, showing that $F$ induces an isomorphism on the $E^1$--pages of the
spectral sequences on both sides amounts to checking that
\[ F': (C_\lft\ten\x_1, \del^0_\lft \ten \mathds{1}) \to (CF^+_\text{left} \otimes
CF^+_\text{right}, \del^0_\lft \ten \mathds{1} + \mathds{1} \ten \del_\rgh) \] is an
isomorphism in homology, where $\del^0_\lft\ten \mathds{1}$ denotes those holomorphic
maps contributing to the differential $\del_{F}$ with $n_w=0$ (Here we have used Lemma
\ref{same multiplicity} to identify $n_w=0$ part of $\del^+$ with $\del^0_\lft\ten
\mathds{1} + \mathds{1} \ten \del_\rgh$). The injectivity of $F'$ in homology follows
from the fact that, by Corollary \ref{Crgh} (see also the proof of Theorem \ref{fibred}
), $\x_1$ does not lie in the image of $\del_\rgh$. Thus, we only need to check that $F'$ is surjective in homology. Suppose that
$\av_1\x_1+\ldots+\av_{4k}\x_{4k}+\bv_1\y_1+\ldots+\bv_{4k}\y_{4k} \in CF^+_\text{left} \otimes
CF^+_\text{right}$ is in the kernel of $\del^0_\lft \ten \mathds{1} + \mathds{1} \ten \del_\rgh$, where
we have chosen the notation so that $\av_i$ and $\bv_i$ are elements in $CF^+_\lft =
C_\lft \ten \Lambda[ \f{Z}_{\geq 0}] $, and $\x_i$ and $\y_i$ are the generators of
$C_\rgh$ as in Theorem \ref{fibred}. Now, because this element is in the kernel of $\del^0_\lft \ten \mathds{1} + \mathds{1} \ten \del_\rgh$, we
have  \begin{eqnarray*}  \del^0_\lft \av_1 = 0 \ \ \ &\text{and}& \ \ \ U\av_1
	+\del^0_\lft \bv_1 = 0 \\
 \del^0_\lft \av_i = 0 \ \ \ &\text{and}& \ \ \ \av_i + \del^0_\lft \bv_i = 0 \ \ \
\text{for} \ \ \ i\neq 1
\end{eqnarray*}

where $U:CF^+_\lft \to CF^+_\lft$ is the usual map in Heegaard Floer theory which
maps $[\av,i] \to [\av,i-1]$. It appears in the above equation because the
disk $D_1$ connecting $\x_1$ to $\y_1$ intersects the base point $z$ with multiplicity
$1$. (Here we also chose an orientation system so that $\del_\rgh \x_i = \y_i$, one can also do the same calculation if $\del_\rgh \x_i = -\y_i$.)

Now, observe that the above equations give 
\begin{eqnarray*} (\del^0_\lft \ten \mathds{1} + \mathds{1} \ten \del_\rgh) (\bv_i\x_i) = - \av_i \x_i +
	\bv_i \y_i \ \ \ \text{for} \ \ \ i\neq 1 \end{eqnarray*}

Therefore, $-\av_i\x_i + \bv_i \y_i$ is in the kernel, but by assumption we
also had $\av_i\x_i + \bv_i \y_i$ in the kernel. This gives us that $2\bv_i
\y_i$ is in the kernel, which in turn, implies that $\del_\lft^0 \bv_i = 0$
(This holds unless we work over a field of characteristic $2$, see below for
that case). Thus, $\av_i =0$ and $\bv_i\y_i$ is in the image of $\del^0_\lft
\ten 1 + 1 \ten \del_\rgh$ (In characteristic $2$, we directly conclude that
$\av_i \x_i+\bv_i \y_i$ is in the image). Therefore, in either case we can
ignore all the terms other than $\av_1\x_1 + \bv_1\y_1$. Furthermore, note that
\[ (\del^0_\lft \ten 1 + 1 \ten \del_\rgh)(U^{-1}\bv_1\x_1)= U^{-1}\del^0_\lft
\bv_1 \x_1 + \bv_1\y_1= -\av_1\x_1 + \bv_1 \y_1 \] Thus, we conclude that
$2\bv_1\y_1$ is in the kernel, which implies that $\del_\lft^0 \bv_1 = 0$
hence, $U\av_1 =0$ and $(\del^0_\lft \ten 1 + 1 \ten
\del_\rgh)(U^{-1}\bv_1\x_1) = \bv_1\y_1$ hence we can ignore the term
$\bv_1\y_1$ and the fact that $U\av_1 =0$ implies that $\av_1\x_1$ is in the
image of $F$ as desired.

This concludes the proof of Proposition \ref{prop1} since a filtered chain map that
induces an isomorphism of $E^1$--pages induces an isomorphism at all pages of the
spectral sequences (see e.g. Theorem 3.5 of \cite{McCleary}), in particular the
$E^\infty$--pages are the groups that we have written in the statement of Proposition
\ref{prop1}. \QED

An immediate corollary that follows from the proof of Proposition \ref{prop1} is that
the $U$-action on $HF^+(Y,f,\gamma_w,\s)$ is trivial for $\s \in
\mathcal{S}(Y|\Sigma_{\text{min}})$. In fact, we have a splitting of
the long exact sequence induced by the $U$-action, which implies the
following relation with the hat-version of Heegaard Floer homology where the
differential counts the holomorphic curves with $n_z=0$ (see \cite{OS}).

\begin{corollary} For $\s \in \mathcal{S}(Y|\Sigma_{\text{min}})$, 
	$\widehat{HF}(Y,f,\gamma_w,\s) \simeq HF^+(Y,f,\gamma_w,\s) \oplus
	HF^+(Y,f,\gamma_w,\s)[1]$ \QED
\end{corollary}

Note that in the case that $g(\Sigma_{\text{min}})=k>1$, there is no perturbation
required thus the above equality holds for the homology groups with $\f{Z}_2$ coefficients. In particular, this implies that $HF^+(Y,\s)$ is algorithmically
computable for $\s \in \mathcal{S}(Y|\Sigma_{\text{min}})$ since there are known algorithms for computing $\widehat{HF}(Y,\s)$.  

Finally, we are ready to state and prove our main result. Over the course of the proof
of the following theorem, we will see why the variant of Heegaard Floer homology that
we denoted by $QFH'(Y,f,\s;\Lambda)$ is well-defined. More precisely, we will see that
the differential that we defined for $QFH'(Y,f,\s;\Lambda)$ squares to zero.

\begin{theorem} \label{isomorphism} $HF^+(Y,f,\gamma_w, \s ) \simeq QFH'(Y,f,\s;\Lambda )$ for $\s\in \mathcal{S}(Y|\Sigma_\text{min})$.
\end{theorem}

\emph{Proof.} Because of Proposition \ref{prop1}, it suffices to prove that \[
H(C_\lft \ten \x_1, \del_\lft \ten \mathds{1}, \gamma_w , \s)  \simeq H(C_\lft, \del ,
\s) \] where the latter group is what we previously called $QFH(Y,f,s)$. Clearly,
we have a one-to-one correspondence between the generators. Next, we will show that
there is an isomorphism of chain complexes. In fact, we will show that the signed
counts of Maslov index $1$ holomorphic curves in
$\pi_2(\x_\lft\ten\x_1,\y_\lft\ten\x_1)$ which contribute to $\del_\lft \ten 1$ and
Maslov index $1$ holomorphic curves in $\pi_2(\x_\lft,\y_\lft)$ that contribute to the
differential $\del$ are equal. First observe that for curves which stay away from the necks at
$\alpha_0$ and $\beta_0$, which are precisely those with $n_w=0$, this one to one
correspondence is clear. (These are the curves that contribute to the differential
$\del^0_\lft \ten 1$ in Proposition \ref{prop1}).

Next, we discuss the curves which have $n_w\neq0$. The correspondence in this case
will be obtained by stretching the necks along $a$ and $b$, which are respectively
parallel pushoffs of $\alpha_0$ and $\beta_0$ to the left of the Heegaard diagram
(into the region $F$).

Let us first describe the holomorphic curves that contribute to $\del_\lft \ten
\mathds{1}$ with $n_w\neq 0$ more precisely. Remember that by definition
$\del_\lft\ten\mathds{1}$ counts those holomorphic curves whose domain lies in $F$, hence they have $n_z=0$. Now, recall that Lemma \ref{same multiplicity} says that by choosing the almost complex structure on $\Sigma \times [0,1] \times \f{R}$ appropriately, one can arrange that the projection to the Heegaard surface is an unbranched cover around the necks $a$ and $b$ (i.e. the holomorphic curve converges to Reeb orbits around $a$ and $b$). Let $A \in \pi_2(\x_\lft \ten \x_1, \y_\lft\ten\x_1)$ be a Maslov index $1$
homology class which is contributing to $\del_\lft \ten \mathds{1}$. By degenerating
the almost complex structure around $a$ and $b$ on $\Sigma$, we get two homology
classes $A_\lft \in \pi_2(\x_\lft,\y_\lft)$ and $A_\rgh \in \pi_2(\x_1,\x_1)$. The
domain of $A_\lft$ lies on $\Sigma_\text{max}$ and it determines a homology class for
the type of holomorphic curves contributing to the differential $\del$. The domain of
$A_\rgh$ has two components $A_\rgh^a$ and $A_\rgh^b$, both supported in disks which
are the domains between $\alpha_0$ and $a$, with $a$ collapsed to a point, and
between $\beta_0$ and $b$ with $b$ collapsed to a point. We claim that the Maslov index of
$A_\lft$ is equal to $1$, and the Maslov indices of each of the components in $A_\rgh$
are equal to $2n_w$. Since the degeneration is along Reeb orbits, we have the formula
\[\text{ind}(A) = \text{ind}(A_\lft) + \text{ind}(A_\rgh^a) +\text{ind}(A_\rgh^b) -
2(N_a + N_b) \] 

where $N_a$ and $N_b$ are the numbers of connected components of the unramified covering in the necks at $a$ and $b$ (clearly $N_a, N_b \in [1,n_w]$).  Therefore, it suffices to see that $\text{ind}(A_\rgh^a)= \text{ind}(A_\rgh^b) =2n_w$.
This follows from the usual formula $\text{ind}(A_\rgh^a)=\langle c_1(\s), A_\rgh^a \rangle =
e(A_\rgh^a) + 2n_\x(A_\rgh^a) = 2n_w$ (since the homology class $A_\rgh^a$ is $n_w$
times the disk with boundary on $\alpha_0$, $e(A_\rgh^a)=n_w$ and $n_\x(A_\rgh^a)=n_w
/2$); similarly for $A_\rgh^b$. We deduce that $\text{ind}(A_\lft) = 1 + 2(N_a+ N_b) -
4 n_w$, which implies that $\text{ind}(A_\lft) = 1$ and the coverings in the
cylindrical necks near $a$ and $b$ are both trivial (in other terms, after
neck-stretching we have $n_w$ distinct cylinders passing through each neck).

Furthermore, we have the evaluation maps : \begin{eqnarray*}
	ev^a_\lft &:& \mathcal{M}(A_\lft) \to \sym^{n_w}(\f{D}) \\  
	ev^a_\rgh &:& \mathcal{M}(A_\rgh^a) \to \sym^{n_w}(\f{D}) \\
	ev^b_\lft &:& \mathcal{M}(A_\lft) \to \sym^{n_w}(\f{D}) \\
	ev^b_\rgh &:& \mathcal{M}(A_\rgh^b) \to \sym^{n_w}(\f{D}) \end{eqnarray*}

given by taking the preimages of the degeneration points of $a$ and $b$ and projecting
to $\f{D}= [0,1]\times\f{R}$. We claim that the moduli space $\mathcal{M}(A)$ can be
identified with the fibre product of moduli spaces $\mathcal{M}(A_\lft) \times_{B}
\mathcal{M}(A_\rgh)$, where $B = \sym^{n_w}(\f{D}) \times \sym^{n_w}(\f{D})$ and the
fibre product is taken with respect to the above evaluation maps.  This is a
consequence of a gluing theorem (see \cite{OSlink} Theorem 5.1 for the proof in a very
closely related situation). 

Finally, we will prove that
$(ev^a_\rgh,ev^b_\rgh):\mathcal{M}(A_\rgh^a)\times \mathcal{M}(A_\rgh^b)\to
\sym^{n_w}(\f{D})\times \sym^{n_w}(\f{D})$ has degree $1$. This implies that, for the
purpose of counting pseudoholomorphic curves, the fibre
product of moduli spaces $\mathcal{M}(A_\lft) \times_{B} \mathcal{M}(A_\rgh)$ can be
identified with $\mathcal{M}(A_\lft)$. Therefore, we can identify the moduli
spaces $\mathcal{M}(A)$ and $\mathcal{M}(A_\lft)$, as required.  

To see that the evaluation maps have degree $1$, we argue as follows: First, we
represent the domain of the strips in $\mathcal{M}(A_\rgh^a)$ by the upper half of the
unit disk so that the upper half circle maps to $\alpha_0$ and the interval
$[-1,1]$ maps to the $\beta$ curve. Also, represent the target disk by the unit disk,
so that $\alpha_0$ corresponds to the unit circle and the $\beta$ arc is represented by
the real positive axis, furthermore the degeneration point of $a$ as used to define
the map $ev^a_\rgh$ is mapped to the origin in this representation. Thus, the moduli
space $\mathcal{M}(A_\rgh^a)$ consists of holomorphic maps from the upper half
disk to the unit disk and $ev^a_\rgh$ records the positions of the $n_w$ zeroes of
these maps. Now, any holomorphic map from the upper half disk to the unit disk can be
reflected ($u(1/\bar{z}):=1/\overline{u(z)}$) to get a holomorphic map from the upper
half-plane to $\f{P}^1$, mapping the real axis to the real positive axis. This can
then be further reflected about the real axis to get holomorphic maps from $\f{P}^1$
to $\f{P}^1$ which are hence rational fractions of degree $2n_w$, with real
coefficients (forced by the invariance under conjugation) and with equivariance under
$z \to 1/\bar{z}$. Now, such holomorphic maps are classified by their zeroes (the
poles are the reflections of the zeroes). In our case, there are $2n_w$ zeroes and
none of these are real, so they are $n_w$ pairs of complex conjugate points. Finally,
we note that $ev^a_\rgh$ maps any such holomorphic map to the positions of its $n_w$
zeroes which lie inside the upper half-disk. Therefore, $ev^a_\rgh :
\mathcal{M}(A_\rgh^a) \to \sym^{n_w}(\f{D})$ is in fact a diffeomorphism. In
particular, it has degree $1$.  \QED

Note that when the minimal genus fibre has genus greater than $1$, there is no
perturbation required since the diagrams that we consider are weakly admissible
in that case. Hence, we get the above isomorphism for the homology groups with
$\f{Z}_2$ coefficients.

\begin{corollary} 
	\label{integral} 
	Suppose that $g(\Sigma_\text{min}) = k >1$, then for
	$\s\in\mathcal{S}(Y|\Sigma_\text{min})$ we have \[ QFH'(Y,f,\s;\f{Z}_2)
	\simeq HF^+(Y,\s). \]
\end{corollary}
        \emph{Proof.} This follows from the above result by letting $t=1$. \QED

In some cases, the quilted Floer homology groups can be calculated easily, the
following special case is an example of this. Given two simple closed curves $\alpha$
and $\beta$ on a surface of genus greater than $1$, let $\iota(\alpha,\beta)$ denote
the geometric intersection number of $\alpha$ and $\beta$, i.e. the number of
transverse intersections of their geodesic representatives for a hyperbolic metric.
\begin{corollary} \label{special case} Suppose that $f$ has only two critical points,
	and let $\alpha , \beta \subset \Sigma_{\text{max}}$ be the vanishing cycles
	for these critical points. Then $\oplus_{\s \in
	\mathcal(S|\Sigma_{min})}HF^+(Y,f,\gamma_w, \s)$ is free
	of rank  $\iota (\alpha,\beta)$.  
\end{corollary} 
\emph{Proof.} When $f$ has only two critical points, $QFH'(Y,f)$ reduces to the Lagrangian
Floer homology of the simple closed curves $\alpha$ and $\beta$ on the surface
$\Sigma_{\text{max}}$. This is easily calculated by representing the free homotopy
classes of simple closed curves $\alpha$ and $\beta$ by geodesics, which ensures that
there are no non-constant holomorphic discs contributing to the differential. In fact,
any holomorphic disk would lift to a holomorphic disk in the universal cover
$\mathbf{H}^2$, which would contradict the fact that there is a unique geodesic
between any two points in $\mathbf{H}^2$.  Therefore, the quilted Floer homology is
freely generated by the number of intersection points of geodesic representatives of
$\alpha$ and $\beta$. \QED

We remark that if $\iota (\alpha , \beta) = 1$ , then the critical values can be
cancelled. Thus for non-fibred manifolds which admit a broken fibration with only $2$
critical points the rank of quilted Floer homology is greater than $1$.

\subsection{An application to sutured Floer homology}

Juh\'asz introduced an extension of Heegaard Floer homology to ``balanced sutured
3--manifolds''. (See \cite{juhasz} for the definition). A connected balanced
sutured manifold $Y$ is a compact oriented 3--manifold with boundary $Y$ and it
can be equipped with a broken fibration $f:Y \to [0,1]$ whose fibers are
surfaces with non-empty boundary and $f^{-1}(0)$ and $f^{-1}(1)$ are
homeomorphic surfaces such that each connected component has exactly one
boundary component (balanced condition). We can always arrange that
$f^{-1}(1/2) = \Sigma$ is the highest genus fibre which is connected and as one
travels from $1/2$ to $0$ one attaches two handles along
$\beta_1,\ldots,\beta_m$ and as one travels from $1/2$ to $1$ one attaches two
handles along $\alpha_1,\ldots,\alpha_m$ which are realized as vanishing cycles
of $f$ on $\Sigma$. The balanced condition translates to the condition that the
set of $\alpha$ curves and respectively the set of $\beta$ curves are linearly
independent in $H_1(\Sigma)$. The sutures $s(\gamma)$ of $Y$ correspond to the
boundary components of $\del \Sigma$ and the annular neighborhoods $A(\gamma)$
of $Y$ are obtained from $s(\gamma)$ by flowing using the gradient flow of $f$
along $[0,1]$ with respect to a metric such that the gradient vector field of
$f$ preserves the boundary of $Y$. 

In \cite{juhasz}, Juh\'asz constructs a variant of Heegaard Floer homology for sutured
$3$--manifolds. This is simply, the Lagrangian Floer homology group
$HF(\sym^{m}(\Sigma), \alpha_1\times\ldots\times\alpha_m,
\beta_1\times\ldots\times\beta_m)$ where the projections of the holomorphic curves
contributing to the differential on $\Sigma_m$ are required to stay away
from the boundary of $\Sigma_m$.

In \cite{KM}, Kronheimer and Mrowka construct an invariant of sutured manifolds using
monopole (resp. instanton) Floer homology, by constructing a closed 3--manifold $Y_n$
and setting the sutured Floer homology of $Y$ by defining it to be a summand of the monopole (resp.
instanton) Floer homology of $Y_n$. The construction of $(Y_n,f_n)$ is by first gluing
$T \times [0,1]$ where $T$ is an oriented connected genus $n\geq 1$ surface with
non-empty boundary, so that $\del
T \times [0,1]$ is glued to the union of annuli $A(\gamma)$, and then identifying the
fibres over $0$ and $1$ by choosing a homeomorphism between them. Note that the
balanced condition implies that $f_n$ has connected fibres. In the monopole (resp.
instanton) setting, Kronheimer and Mrowka define the sutured monopole (resp.
instanton) Floer homology of $Y$ to be $\bigoplus_{s\in \mathcal{S}(Y| \Sigma_{\text{min}})}HM(Y_n, \s)$ and prove that this is an invariant
of the sutured manifold $Y$ (in particular, it is also independent of the genus $n$ of
$T$ and the homeomorphism chosen in identifying fibres over $0$ and $1$). It was raised in \cite{KM} as a question, whether one can recover Juh\'asz's
definition of sutured Floer homology from the construction given above applied in the
setting of Heegaard Floer homology. In the next theorem, we give an affirmative answer
to this. 

\begin{theorem} For $n\geq 1$,
	\[  SFH(Y) \simeq \bigoplus_{\s\in \mathcal{S}(Y_n|
	\Sigma_{\text{min}})} HF^+(Y_n, \s) \] 
\end{theorem} 
	
Note that this theorem in particular implies that the group on the right hand side is
independent of $n$ and the chosen surface homeomorphism in the construction of $Y_n$.
As usual, in the case that the lowest genus fibre of $f_1$ has genus $1$, one needs to
use coefficients in $\Lambda$. 

\emph{Proof.} Theorem $\ref{isomorphism}$ applied to $(Y_n,f_n)$ yields that
$\bigoplus_{\s\in\mathcal{S}(Y_n| \Sigma_{\text{min}})} HF^+(Y_n,\s)  =
QFH'(Y_n,f_n)$.  Therefore, the proof will follow once we establish that
$SFH(Y,f) \simeq QFH'(Y_n,f_n)$. This in turn relies on a simple observation
about the Heegaard diagrams used in the definition of these groups, namely let
us denote the maximal genus fibre of $f$ by $\Sigma$, and the maximal genus
fibre of $f_n$ by $\Sigma \cup T$. Now, if an admissible sutured Heegaard
diagram of $(Y,f)$ is given by $(\Sigma, \alpha_1,\ldots,\alpha_m,
\beta_1,\ldots,\beta_m)$, then the Heegaard diagram for calculating
$QFH'(Y_n,f_n)$ is given by $(\Sigma\cup T, \alpha_1,\ldots,\alpha_m,
\beta_1,\ldots,\beta_m)$.  Note that there is no $\alpha$ or $\beta$ curve
entering $T$. Thus, the proof will be complete once we prove that holomorphic
curves contributing to the differential of $QFH'(Y_n,f_n)$ do not enter to the
region including $T$. Note that because of the admissibility condition of the
sutured Heegaard diagram of $(Y,f)$ we can use an almost complex structure
which is vertical in a neighborhood of $\Sigma \times [0,1] \times \f{R}$ (by
vertical, we mean that the fibres of the projection $\Sigma \times [0,1] \times
\f{R} \to [0,1] \times \f{R}$ are holomorphic) so that the holomorphic curves
contributing to the differential of sutured Floer homology appear as
holomorphic curves contributing to the differential of $QFH'(Y_n, f_n)$.  On
the other hand, we use a non-vertical almost complex structure as in Section
\ref{variant}, along $T\times [0,1] \times \f{R}$ away from the boundary of
$T$. Now, let $u: (S,\del S) \to  (\Sigma\cup T) \times [0,1] \times \f{R}$ be
a holomorphic map contributing to the differential of $QFH'(Y_n,f_n)$. We would
like to show that the image of the projection of $u$ to the Heegaard surface
does not hit $T$. This follows from a degeneration argument. Namely, suppose
that the image of the projection of $u$ does hit $T$, then we can degenerate
along Reeb orbits corresponding to the attaching region of $T$ to $\Sigma$,
this would on one side give a holomorphic map $u_T : \tilde{S} \to \tilde{T}
\times [0,1] \times \f{R}$ where $\tilde{T}$ is the closed surface obtained by
shrinking each boundary component of $T$ to a point and $\tilde{S}$ is a part
of the domain of the degenerated map. It follows for example from Corollary 4.3
in \cite{lipshitz} that the index of $u_T$ is equal to the Euler measure of the
projection onto $\tilde{T}$. (Strictly speaking Corollary 4.3 in
\cite{lipshitz} is proved to hold for holomorphic curves with corners. Here we
apply it in the degenerate case where there are no corners. This can be
justified as follows: degeneration of $u$ along the attaching region results in
$u$ degenerating into two pieces $u_\Sigma$ and $u_T$. Then Corolarry 4.3 in
\cite{lipshitz} can be applied to both $u$ and $u_\Sigma$ and a short
calculation using the additivity of the index then implies that index of $u_T$
is the Euler measure of $\tilde{T}$). Furthermore, in this case (since there
are no corners) the Euler measure is simply $\chi(\tilde{T})$ times the
multiplicity of the domain supported in the whole of $\tilde{T})$. The
multiplicity is positive by holomorphicity of $u_T$ and $\chi(\tilde{T})$ is
negative by assumption.  So, we conclude that the index is negative.
Furthermore, we have restricted to the class of almost complex structures so that the fibre of
the projection $(\Sigma \cup T) \times [0,1] \times \f{R} \to (\Sigma \cup T)
\times [0,1]$ is not a holomorphic surface, this ensures that one can choose an almost complex structure giving transversality as in
Proposition 3.7 of $\cite{lipshitz}$. This yields the desired contradiction
(note that in the case that $T$ has genus $1$, we still obtain a contradiction
since we get a negative dimension for the transversely cut moduli space after
quotienting by the $\f{R}$ action). \QED 

\section{Isomorphism between $QFH(Y,f;\Lambda)$ and $QFH'(Y,f;\Lambda)$ }
\label{section4}

In this section, we relate $QFH'(Y,f)$ defined as a variant of Heegaard Floer homology
as in Section \ref{variant} with the original definition in terms of holomorphic
quilts given in the introduction, which we called $QFH(Y,f)$ (see below for a detailed definition). More precisely, we show
that these two groups are isomorphic whenever they are defined. 

There are two main ingredients in this isomorphism. The first one is a general
result in the theory of holomorphic quilts which proves an isomorphism between
quilted Floer homology groups under transverse and embedded compositions of
Lagrangians (see below for definitions). This result is originally due to
Wehrheim and Woodward \cite{katrin}, which was proved in the positively
monotone setting. However, here we are situated in the (strongly) negative
setting in which case the arguments of Wehrheim and Woodward are no longer
valid. To resolve this issue, in \cite{LL} we gave a new proof of Wehrheim and
Woodward's theorem which applies in the current situation. The new proof applies
under both positive and (strongly) negative monotonicity assumptions. 

The second main ingredient in the proof of the isomorphism is a detailed study
of the Lagrangian correspondences that are involved in the definition of
$QFH(Y,f)$. In the next section, we give a detailed definition of $QFH(Y,f, \s )$  for $\s \in
\mathcal{S}(Y|\Sigma_k)$ and $g<2k$. In particular, we give a detailed description of monotonicity which is required to have a rigorous definition over $\f{Z}_2$ when $k>1$. Proving the isomorphism involves showing that various
compositions of these Lagrangians correspondences are Hamiltonian isotopic to
product tori, that appear in the definition of $QFH'(Y,f)$.  This part of the
proof has appeared in author's thesis \cite{lekiliThesis}, and it can also be
found in the upcoming work \cite{lekiliPerutz} where a more general set-up is
developed.

Finally, we remark that all the theorems are stated for Floer homology groups
over the universal Novikov ring $\Lambda$, but as before, in the case where the
lowest genus fibre has genus greater than $1$, one can take coefficients to be
in $\f{Z}_2$. 

\subsection{Definition of $QFH(Y,f,\s)$}
\label{qfhdef}

We now give a detailed definition of $QFH(Y,f,\s)$ for $\s \in
\mathcal{S}(Y|\Sigma_\text{max})$ when $g(\Sigma_\text{max}) <
2g(\Sigma_\text{min})$.  Recall that we start with a broken fibration $f:Y \to
S^1$ with connected fibres and with a distinguished maximal genus fibre $\Sigma_g = f^{-1}(-1)$
and a minimal genus fibre $\Sigma_k = f^{-1}(1)$. There are $g-k$ critical values $p_1, \ldots,
p_{g-k}$ on the northern semi-circle in clockwise order and $g-k$ critical
values $q_1,\ldots , q_{g-k}$ on the southern semi-circle in counter-clockwise
order. For a critical value $p$ fix two nearby points $p^+,p^-$ on $S^1$ such
that the genus $f^{-1}(p^+)$ is greater than the genus of $f^{-1}(p^-)$, that
is, $p^+$ is to the left of $p^-$. Furthermore, we arrange that $p^+_i =
p^-_{i+1}$ and that $p_1^+ = q_1^+ = -1$ and $p_{g-k}^- \neq  q_{g-k}^- $. 

Next choose a Riemannian metric $g$ on $Y$, we then have embedded curves
$\alpha_i \subset f^{-1}(p_i^+)$ and $\beta_i \subset f^{-1}(q_i^+)$ cut out by
the unstable manifolds of $p_i$ and $q_i$. By abuse of notation, we also denote by
$\alpha_i, \beta_i \subset \Sigma_g = f^{-1}(-1)$ the embedded curves cut out
by the intersection of the unstable manifolds of $p_i$ and $q_i$ with
$\Sigma_g$. 

Finally, choose an area form $\xi_i$ and a compatible complex structure $j_i$
on each $f^{-1}(p_i^+)$ , $f^{-1}(q_i^+)$ and $f^{-1}(1)$. Note that the
gradient flow gives an identification of $f^{-1}(p_i^+)$ and $f^{-1}(p_i^-)$
(resp. for $f^{-1}(q_i^+)$ and $f^{-1}(q_i^-)$) outside of their intersection
with the stable and unstable manifolds associated with $p_i$ (resp. $q_i$) and we
ask that this identification is a complex isomorphism when the fibres are
equipped with the complex structures. Note that $f^{-1}(p_{g-k}^-)$ and
$f^{-1}(q_{g-k}^-)$ are diffeomorphic surfaces and we can and do in fact arrange them to be
symplectomorphic, however we cannot in general demand that they are isomorphic
complex surfaces as this will put severe restrictions on the monodromy.

Let $B = S^1 - \left( \left( \bigcup_{i=1}^{g-k} ( (p_i^+,p_i^-) \cup (q_i^+,q_i^-) \right) \cup (p_{g-k}^- , q_{g-k}^-) \right) $ be
the finite set of $2(g-k)+1$ points, one between every consecutive pair of
critical points with the exception of $p_{g-k}^-$ and $q_{g-k}^-$ which lie between the consecutive critical points $p_{g-k}$ and $q_{g-k}$.
Except $p_{g-k}^-$ and $q_{g-k}^-$, any two consecutive points (as a subset of $S^1)$ gives a
quintuple $(\Sigma, j, C, \overline{\Sigma}, \overline{j})$ where $(\Sigma,j)$
and $(\overline{\Sigma}, \overline{j})$ are connected Riemann surfaces, $C$ is
an embedded non-separating curve on $\Sigma$ and there is a canonical
diffeomorphism from $\Sigma_C$ (the result of surgery on $C$) to 
$\overline{\Sigma}$.

For $b \in B$, let $F_b = f^{-1}(b)$ be the fibre of $f$ equipped with a
complex structure and a compatible area form $\xi_b$ as above. We can then
consider $\sym^n(F_b)$ as a complex manifold for any $n$. There are two
distinguished classes in $H^2(\sym^n(F_b) )$ which span the invariant subspace
of the second cohomology group under the action of the mapping class group of
$F_b$. These are $\eta$, Poincar\'e dual to $\{pt\} \times \sym^{n-1}(F_b)$ and
$\theta$, which can be concisely described using the fact that
$c_1(T\sym^n(F_b)) =(n+1-g_b) \eta - \theta$, where $g_b$ is the genus of
$F_b$. 

In \cite{handleslide}, Perutz constructs K\"ahler forms $\omega_{F_b}$ on
$\sym^n(F_b)$ with the property that $\omega_{F_b}$
agrees with $\sym^n{\xi_b}$ outside of a neighborhood
of the diagonal and $[\omega_{F_b}] = \eta + \lambda
\theta$ for any sufficiently small fixed real parameter $\lambda
>0$, and which tames $\sym^n(j)$ on all of $\sym^n(F_b)$. 

We are now ready to state the fundamental construction of Perutz:

\begin{theorem} \label{per} (Perutz \cite{LM1}) Starting from a quintuple $(\Sigma, j, C,
	\overline{\Sigma}, \overline{j})$ as above, one can construct a
	Lagrangian correspondence $L_C \subset \sym^{n-1}(\overline{\Sigma} )
	\times \sym^n (\Sigma) , -\omega_{\overline{\Sigma}} \oplus
	\omega_\Sigma )$ canonically up to Hamiltonian isotopy.  \end{theorem}

	We note here two topological properties of
	$L_C$ from \cite{LM1} : First, there are maps \[
	\sym^n(\Sigma) \xleftarrow{i} L_C
	\xrightarrow{\pi} \sym^{n-1}(\overline{\Sigma})
	\]	such that $i$ is a codimension $1$
	embedding and $\pi$ is a trivial $S^1$
	fibration. Second, note that for $n>1$
	$\pi_1(\sym^n(\Sigma)) = H_1(\Sigma)$ and the
	homology class of $L_C$ in $\sym^n(\Sigma)
	\times \sym^{n-1}(\overline{\Sigma})$ is given
	by $C \times \sym^{n-1}(\overline{\Sigma})$.

	For $\s$ a $\text{spin}^c$ structure on $Y$, let us define the integer $n_b$ by the
	formula $\langle c_1(\s) , F_b \rangle = 2n_b + \chi(F_b)$.

Perutz's construction applied to our set-up gives a sequence of Lagrangian
correspondences between $\sym^{n_{p_{g-k}^-}}(F_{p_{g-k}^-})$ and
$\sym^{n_{q_{g-k}^-}}(F_{q_{g-k}^-})$. Furthermore, these latter two symplectic
manifolds are canonically identified by a symplectomorphism induced by the
symplectomorphism of $F_{p_{g-k}^-}$ and $F_{q_{g-k}^-}$. (The fact that the
complex structures on $F_{p_{g-k}^-}$ and $F_{q_{g-k}^-}$ are compatible with
the same symplectic structure provides a path of complex structures
interpolating between the given two complex structures on the underlying
surface, this in turn gives a tautological K\"ahler isomorphism between
$\sym^{n_{p_{g-k}^-}}(F_{p_{g-k}^-})$ and
$\sym^{n_{q_{g-k}^-}}(F_{q_{g-k}^-})$. See for example \cite{usher} for
further details of this identification where a definition of the Floer
homology group that we are discussing here was given for the special of fibred
3--manifolds). In this paper, we are only concerned with the spin$^c$
structures $\s \in \mathcal{S}(Y|\Sigma_k)$. Thus, for $\s \in
\mathcal{S}(Y|\Sigma_k)$, we have $2n_b = 2-2k-\chi(F_b)$. Hence, we have
$n_{p_{g-k}^-} = n_{q_{g-k}^-} = 0$ and so the above identification is trivially the identity map. 

In any case, with the above identification in mind,  we obtain a cyclic set of
Lagrangian correspondences $L_{\alpha_1}, \ldots, L_{\alpha_{g-k}}$ and
$L_{\beta_1}, \ldots, L_{\beta_{g-k}}$ between the cyclic set of symplectic
manifolds $\{ \sym^{n_b} (F_b) \}_{b\in B}$.

Starting from such data, we define the quilted Floer homology of $(Y,f)$ as the
quilted Floer homology of this cyclic set of Lagrangians as developed by
\cite{katrin} (see also \cite{gysin}) : \[ QFH(Y,f) := HF(L_{\alpha_1}, \ldots,
L_{\alpha_{g-k}}, L_{\beta_{g-k}},\ldots, L_{\beta_1}) \] 

Let us recall the basic definition of quilted Floer
homology of a cyclic set of Lagrangians, as we will
need some of this notation in order to prove that our
quilted Floer homology is well-defined. Let us choose
a cyclic ordering of the set $B$, write $b_i$ for the $i^{th}$ element (where the indices are always considered in $mod\ 2(g-k)$) and write
$\underline{L} = (L_1,L_2,\ldots L_{2(g-k)}) =
(L_{\alpha_1}, \ldots, L_{\alpha_{g-k}},
L_{\beta_{g-k}} , \ldots L_{\beta_1} )$ so that $L_i \subset \sym^{n_{b_i}}(F_{b_i}) \times (\sym^{n_{b_{i+1}}}(F_{b_{i+1}}))$

The quilted
Floer chain complex $CF(\underline{L})$ is freely
generated over the base ring by the generalized
intersection points: \[ I(\underline{L}) = \{
\underline{x} = (x_1,\ldots x_{2(g-k)} | (x_{2(g-k)}, x_1) \in
L_1, (x_1,x_2) \in L_2, \ldots (x_{2(g-k)-1}, x_{2(g-k)}) \in L_{2(g-k)} \} \] This set can be arranged to be finite by
requiring the curves $\alpha_i$ and $\beta_j$ intersect
at finitely many points, for all $i, j$. (This will be
made clear below). Next consider the path space \[
\mathcal{P}(\underline{L}) = \{ (\gamma_1,\ldots, \gamma_{2(g-k)})
| \gamma_i:[0,1] \to \sym^{n_{b_i}}(F_{b_i}) , (\gamma_i(1), \gamma_{i+1}(0)) \in L_i \} \]

Note that all the symplectic manifolds $\sym^{n_{b_i}}(F_{b_i})$ are equipped with symplectic forms $\omega_i$ in the class $\eta+\lambda \theta$ and we choose compatible almost complex structures $J_i$. The Floer differential is then obtained in the usual way by counting moduli space of finite energy quilted holomorphic strips connecting generalized intersection points $\underline{x}$ and $\underline{y}$ as follows: 
\begin{eqnarray*}
	\mathcal{M}(\underbar{x},\underbar{y}) = \{ u_i : \f{R} \times [0,1] \to \sym^{n_{b_i}}(F_{b_i}) |
	\ \delbar_{J_i} u_i = 0 ,\\
	  E(u_i) = \int u_i^*\omega_i < \infty \\ \text{lim}_{s\to-\infty}
	  u_i(s,\cdot) = x_i , \text{\ lim}_{s\to +\infty} u_i(s,\cdot) = y_i \\
	   (u_i(s, 1), u_{i+1}(s,0)) \in L_i\  \text{for all}\ i=1,\ldots 2(g-k) \} /
	   \f{R}
\end{eqnarray*}
For people shy of holomorphic quilts, one can
alternatively think of the latter group as a Floer
homology group of two Lagrangians $\mathbb{L}_1= L_1
\times L_3 \times \ldots \times L_{2(g-k)-1}$ and
$\mathbb{L}_2 = L_2 \times L_4 \times \ldots \times
L_{2(g-k)}$in the product symplectic manifold
$\mathbb{M} = \prod_{b\in B} \sym^{n_b} (F_{b})$.

The quilted Floer homology group is defined under monotonicity assumptions and we have to show that our set-up falls into (strongly) negative monotone case (compare \cite{gysin}, Definition 1.8). We address these technicalities now:

\emph{Transversality and avoiding bubbles:} This follows from
standard arguments in Floer theory, see for example,
\cite{seidel} Lemma 2.4. In the strongly negative monotone
case, in addition to transversality for moduli space of
Floer trajectories, transversality for the moduli space of bubbles can be achieved by an identical argument (see
\cite{gysin} Lemma 3.5 ). We insist on using a path of almost complex structures which are of the form $\sym^{n_b}(j_s)$ near the diagonal as in \cite{OS} on each component $\sym^{n_b}(F_b)$ of $\f{M}$.
Though, we warn the reader that one cannot necessarily achieve transversality by considering complex structures $J_s$ on $\mathbb{M}$ which is a product of generic complex structures on each factor (see \cite{WWcorrect} for a discussion of this issue). Therefore, outside of the neighborhood of the diagonal one uses generic complex structure on $\mathbb{M}$ (i.e., not of split-type.). 

To avoid disk or sphere bubbles, we pick a generic complex structure from our class of almost complex structure described above, which achieves transversality for Floer trajectories as well as disk and sphere bubbles. Then under the assumption $g<2k$ (this is the strong negativity assumption), one can calculate the dimension of the moduli space of disk and sphere bubbles and get a negative number. In view of transversality, this proves that there are no non-constant bubbles. The calculation of the dimension of the moduli space of disk and sphere bubbles follows from Section 4 of \cite{LM2}. We spell out this calculation here for completeness:

Note that the symplectic manifolds that we are dealing with are $M=
\sym^n(\Sigma)$ where $n = g(\Sigma)-k$ and $g(\Sigma)$ takes values between
$g(\Sigma_{\text{max}})=g$ and $g(\Sigma_{\text{min}})=k$. We equipped $M$ with a
symplectic form in the class $\eta+\lambda \theta$ where $\lambda>0$ is a fixed
parameter that is determined by the monotonicity condition as follows:

The monotonicity constant $\tau$ is determined by the equation:
\[ [\eta + \lambda \theta] = \tau [c_1(\sym^n(\Sigma))] = \tau [(n+1-g)\eta -\theta] \]

Therefore, $\tau= \frac{1}{n+1-g} <0$ is the fixed monotonicity constant which is the
same for any of the symplectic manifolds we consider since $n-g=
-g(\Sigma_{\text{min}})$.

Now, Perutz
calculates in Section 4 of \cite{LM2} that for $n>1$ the Hurewicz map $\pi_2(\sym^n(\Sigma))
\to H_2(\sym^n(\Sigma))$ has rank $1$ and generated by a class $h$ which satisfies
$\eta( h)=1$ and $\theta(h)=0$. On the other hand, $c_1(\sym^{n}(\Sigma))= (n+1
-g(\Sigma)) \eta - \theta $. Therefore, any simple holomorphic sphere would have
$[u]=h$ and its index would be : \[ 2(\langle c_1(\sym^{n}(\Sigma) , h \rangle + n -3)
= 4n-2g(\Sigma)-4 = 2g(\Sigma) -4k -4 \] 

The assumption $g<2k$ now implies that this quantity is strictly less than $-4$, which suffices for our purpose. (For $n=1$, we can't have any holomorphic spheres since $\pi_2(\Sigma)=0$) 

Similarly, for a disk bubble we need to verify the assumptions for our Lagrangian $L_C
\subset \sym^{n}(\Sigma) \times \sym^{n-1}(\overline{\Sigma})$, where as before
$n=g(\Sigma)-k=g(\overline{\Sigma})+1-k$. In light of the fact that Perutz proves in Lemma 3.18
of \cite{LM1} that any disk in $\pi_2(\sym^{n}(\Sigma) \times \sym^{n-1}(\overline{\Sigma}),
L_C)$ lifts to a sphere it follows that \[ \mu_{L_C}([u]) = 2(\langle
c_1(\sym^{n}(\Sigma)\times \sym^{n-1}(\overline{\Sigma})) , [u] \rangle)  \] Now,
the positive area disks $u$ for which the value $\mu_{L_C}([u])$ is maximal, have
index given by \[ 2(n+1-g(\Sigma)) + (2n-1) -3  = 2g(\Sigma)-4k-2 \] 
Again, the assumption $g<2k$ ensures that this value is strictly less than $-2$, which
guarantees that the non-existence of disk bubbles in the relevant moduli space of index 0,1 and 2.

\emph{Monotonicity (admissibility):} We have seen that the
symplectic manifolds $\sym^b(F_b)$ are negatively
monotone with the same monotonicity constant $\tau = \frac{1}{2-2k}<0$
as $2c_1(T\sym^b(F_b)) = 2(n_b-g_b+1) \eta - 2\theta$ and
$\omega_{F_b} = \eta + \lambda\theta $, and that each Lagrangian $L_{\alpha_i}$ and $L_{\beta_i}$, hence
their product $\mathbb{L}_1$ and $\mathbb{L}_2$ are
monotone with the same constant $\tau$. Note that the value of $\lambda>0$ is irrelevant since $\theta$ vanishes on spherical classes and any disk with boundary on $\mathbb{L}_i$ come from a spherical class [\cite{LM2}, Section 4]. 

The only missing ingredient is that monotonicity for
the pair $(\mathbb{L}_1, \mathbb{L}_2)$. That is to
say, we need to show that index and the area functions
on $\pi_1(\mathcal{P}(\underline{L}))$ are proportional
with the monotonicity constant $\tau$ as above. This is needed in two places, first we need to have an a priori energy bound for low index moduli spaces in order to appeal to Gromov compactness theorem to say that at the boundary of moduli spaces we either get broken trajectories or bubbled configurations. The bubbled configurations are then eliminated using the strong negativity assumptions. The second place where monotonicity is needed is to show that there are only finitely many homotopy classes of disks between given two intersection points, hence the Floer differential can be defined.

In fact, monotonicity for the pair does not always hold and depends on
the relative position of the curves $\alpha_i$ and $\beta_i$ and our task is to
show that it can be ensured whenever the diagram $(\Sigma_g, \alpha_1,\ldots,
\alpha_{g-k}, \beta_1,\ldots,\beta_{g-k})$ is an admissible diagram in the
sense defined below (see also Section \ref{section2}). 

Our strategy will be to show that admissibility implies monotonicity for the
Heegaard tori $(\mathbb{T}_\alpha,\mathbb{T}_\beta)$ in $\sym^{g-k}(\Sigma_g)$
for a symplectic form $\omega_\xi$ in the class $\eta$ which tames the same complex structure as $\omega$ 
and deduce from that the required monotonicity properties for the pair $(\mathbb{L}_1,\mathbb{L}_2)$ 

We will postpone this until we put our Lagrangians $L_{\alpha_i}$ and $L_{\beta_i}$ in a nice position (by a Hamiltonian isotopy) so as to relate them to Heegaard tori $\alpha_1\times\ldots\times\alpha_{g-k}$ and $\beta_1\times\ldots\times\beta_{g-k}$ in $\sym^{g-k}(\Sigma_g)$. 

\subsection{Heegaard tori as composition of Lagrangian correspondences}

Recall that given two Lagrangian correspondences, $L_1 \subset X\times Y$ and $L_2
\subset Y \times Z$ such that $L_1 \times L_2$ is transverse to the diagonal in $Y$, the composition $L_1 \circ L_2$ is a Lagrangian correspondence
in $X \times Z$ given by the union of tuples $(x,z)$ such that there exists a $y \in
Y$ with the property that $(x,y) \in L_1$ and $(y,z) \in L_2$. 

Now, for the class of almost complex structures $j$ that are sufficiently stretched
along the vanishing cycles of $f$ near its critical points, we have the following important technical lemma about these correspondences which was conjectured by Perutz in \cite{LM1}:

\begin{lemma} \label{lemham} For $g>k$, $
	L_{\alpha_1} \circ \ldots \circ
	L_{\alpha_{g-k}}$ and $L_{\beta_1} \circ \ldots
	\circ L_{\beta_{g-k}}$ are respectively
	Hamiltonian isotopic to $\alpha_1 \times \ldots
	\times \alpha_{g-k}$ and $ \beta_{1} \times
	\ldots \times \beta_{g-k}$ in
	$\sym^{g-k}(\Sigma) $ equipped with a K\"ahler
	form $\omega$ which lies in the cohomology
	class $\eta + \lambda \theta$ with $\lambda
	>0$. \QED \end{lemma}  

The proof of this lemma can be found in
\cite{lekiliPerutz} and \cite{lekiliThesis}. The proof
is obtained by carrying out the construction of
Lagrangian correspondences as a family of
degenerations. As the required technical set-up is
developed extensively in \cite{lekiliPerutz}, for the
sake of brevity we choose to omit it from here.

\emph{Back to periods and admissibility:} Lemma
\ref{lemham} is accomplished by carrying out
the construction of $L_{\alpha_i}$ simultaneously which
enables us to show that for a careful choice of
degeneration one in fact has the exact equality:
\begin{eqnarray*}  L_{\alpha_1} \circ \ldots \circ
	L_{\alpha_{g-k}} &=& \alpha_1\times \ldots
	\times \alpha_{g-k} \\ L_{\beta_1} \circ \ldots
	\circ L_{\beta_{g-k}} &=& \beta_1 \times \ldots
	\beta_{g-k} \end{eqnarray*}

	From now on, we will work with this situation
	and prove that the quilted Floer homology is well-defined in
	an admissible situation. We can then appeal to Hamiltonian isotopy provided by Lemma \ref{lemham} to conclude that
	quilted Floer homology will be well-defined
	even in the cases where the above equalities
	might not hold. (Monotonicity still holds after Hamiltonian isotopy).

Let us denote the Heegaard tori in
$\sym^{g-k}(\Sigma_g)$ by $\mathbb{T}_\alpha$ and
$\mathbb{T}_\beta$ and the path space connecting these
by $\Omega(\mathbb{T}_\alpha, \mathbb{T}_\beta)$. Note
that the above equalities imply that there is bijection
between the $CF(L_{\alpha_1},
\ldots, L_{\alpha_{g-k}}, L_{\beta_{g-k}},\ldots,
L_{\beta_1})$ (generated by $\mathbb{L}_1 \cap \mathbb{L}_2$)  and $CF(L_{\alpha_1} \circ \ldots \circ
L_{\alpha_{g-k}}, L_{\beta_1} \circ \ldots \circ
L_{\beta_{g-k}} )$ (generated by $\mathbb{T}_\alpha \cap
\mathbb{T}_\beta$). Recall that we denoted the path
space that is used in the definition of quilted Floer
homology by $\mathcal{P}(\underline{L})$. This is canonically identified with the path space $\Omega(\mathbb{L}_1, \mathbb{L}_2)$ for the Lagrangians $\mathbb{L}_1$ and $\mathbb{L}_2$. The main topological lemma about these path spaces is the following:

\begin{lemma} \label{lem1} There exists an inclusion map $\iota :
	\Omega(\mathbb{T}_\alpha, \mathbb{T}_\beta) \to
	\Omega(\mathbb{L}_1, \mathbb{L}_2) =
	\mathcal{P}(\underline{L})$, which induces a bijection between $\mathbb{T}_\alpha \cap
\mathbb{T}_\beta$ and $\mathbb{L}_1 \cap \mathbb{L}_2$ and an isomorphism $\pi_1(\Omega(\mathbb{T}_\alpha, \mathbb{T}_\beta), \mathbf{x}) \to
	\pi_1(\Omega(\mathbb{L}_1, \mathbb{L}_2), \iota({\mathbf{x}}))$ for any $\mathbf{x} \in \mathbb{T}_\alpha \cap
\mathbb{T}_\beta$.
\end{lemma} 
	
\emph{Proof.} Let $\gamma:[0,1] \to \sym^{g-k}(\Sigma_g)$ be
	a path in $\Omega(\mathbb{T}_\alpha,
	\mathbb{T}_\beta)$. In particular, $\gamma(0)
	\in \mathbb{T}_\alpha = L_{\alpha_1} \circ
	\ldots \circ L_{\alpha_{g-k}}$. As our
	Lagrangian correspondences are circle bundles
	$\pi_i : L_{\alpha_i} \to
	\sym^{g-k-i}(F_{b_i})$, $\gamma(0)$ determines
	a tuple $(\pi_1(\gamma(0)),
	\pi_2\circ\pi_1(\gamma(0)),\ldots,
	\pi_{g-k}\circ\ldots\circ \pi_1(\gamma(0)))$,
	call these $(x_2,\ldots x_{g-k})$. Similarly,
	$\gamma(1) \in \mathbb{T}_\beta$ determines a
	tuple $(x_{2(g-k)}, x_{2(g-k)-1},\ldots
	x_{{g-k}+1})$. The map $\iota :
	\Omega(\mathbb{T}_\alpha, \mathbb{T}_\beta) \to
	\mathcal{P}(\underline{L})$ is given by: \[
	\gamma \to (\gamma, x_2,x_3,\ldots, x_{2(g-k)})
	\]

In other words, all but the first component are constant paths which are in turn determined by the first component. 

Now, the statement about the bijection between $\mathbb{T}_\alpha
\cap \mathbb{T}_\beta$ and $\mathbb{L}_1 \cap \mathbb{L}_2$ follows from the definition of the map $\iota$. Let us fix an intersection point $\mathbf{x} \in \mathbb{T}_\alpha \cap \mathbb{T}_\beta$. We also write $\mathbf{x}$ for the corresponding intersection point in $\mathbb{L}_1 \cap \mathbb{L}_2$.

Now, consider the path component of $\Omega(\mathbb{T}_\alpha, \mathbb{T}_\beta)$ containing $\mathbf{x}$, one has evaluation maps on both sides which induces the Serre fibration: 
\[ \Omega_{\mathbf{x}} (\sym^{g-k}(\Sigma_g)) \to \Omega(\mathbb{T}_\alpha, \mathbb{T}_\beta) \to \mathbb{T}_\alpha \times \mathbb{T}_\beta \]

There is a similar Serre fibration for $\Omega(\mathbb{L}_1, \mathbb{L}_2)$. The map $\iota$ we defined above now gives a map between the two Serre fibrations:

\begin{displaymath} \xymatrix { \Omega_{\mathbf{x}}
	(\sym^{g-k}(\Sigma_g)) \ar[r]\ar[d]^{\iota}
	&\Omega(\mathbb{T}_\alpha, \mathbb{T}_\beta)
	\ar[r]\ar[d]^{\iota} & \mathbb{T}_\alpha \times
	\mathbb{T}_\beta \ar[d]^{\iota} \\ \Omega_{\mathbf{x}}
	(\mathbb{M}) \ar[r] & \Omega(\mathbb{L}_1,
	\mathbb{L}_2) \ar[r] & \mathbb{L}_1 \times
	\mathbb{L}_2 &} \end{displaymath}

Note that the leftmost arrow can also be seen as induced from the based inclusion of $\sym^{g-k}(\Sigma_g)$ to $\mathbb{M}$. We next consider the long exact sequences induced by these Serre fibrations. We have the following commutative diagram:

\begin{displaymath} \xymatrix {  0 \ar[r]\ar[d] &
	\pi_2(\sym^{g-k}(\Sigma_g))
	\ar[r]\ar[d]^{\iota}
	&\pi_1\Omega(\mathbb{T}_\alpha,
	\mathbb{T}_\beta) \ar[r]\ar[d]^{\iota} &
	\pi_1(\mathbb{T}_\alpha \times
	\mathbb{T}_\beta) \ar[r] \ar[d]^{\iota} &
	\pi_1(\sym^{g-k}(\Sigma)) \ar[d]^{\iota} \\ \pi_2(\mathbb{L}_1
	\times \mathbb{L}_2) \ar[r]^{i} & \pi_2
	(\mathbb{M}) \ar[r] & \pi_1\Omega(\mathbb{L}_1,
	\mathbb{L}_2) \ar[r] & \pi_1(\mathbb{L}_1
	\times \mathbb{L}_2) \ar[r]^{p} & \pi_1(\mathbb{M})
	} \end{displaymath}

From which, one can obtain the following set of short-exact sequences:

\begin{displaymath} \xymatrix {  0 \ar[r]\ar[d] &
	\f{Z}
	\ar[r]\ar[d]^{\iota}
	&\pi_1\Omega(\mathbb{T}_\alpha,
	\mathbb{T}_\beta) \ar[r]\ar[d]^{\iota} &
	\frac{H_1(\Sigma_g)}{[\alpha_1],\ldots,[\alpha_{g-k}], [\beta_1],\ldots, [\beta_{g-k}]} \ar[r] \ar[d]^{\iota} &
	0 \ar[d] \\ 0 \ar[r] & \text{coker } i \ar[r] & \pi_1\Omega(\mathbb{L}_1,
	\mathbb{L}_2) \ar[r] & \ker{p} \ar[r] & 0
	} \end{displaymath}

	From the topological properties of
	$L_{\alpha_i}$ and $L_{\beta_j}$ mentioned
	after Theorem \ref{per} it follows that all
	except the middle arrow is an isomorphism.
	We appeal to five-lemma to conclude that
	the middle arrow is also an isomorphism, as
	desired. \QED

Recall that our goal is to show that the index and the area are proportional on
$ \pi_1(\Omega(\mathbb{L}_1, \mathbb{L}_2), \mathbf{x})$ for the Lagrangian
intersection problem $(\mathbb{M}; \mathbb{L}_1, \mathbb{L}_2)$. We next show
that in view of the above lemma, it suffices to show that the index and the are
are proportional on $\pi_1(\Omega(\mathbb{T}_\alpha, \mathbb{T}_\beta), \mathbf{x}) $ for the Lagrangian intersection problem $(\sym^{g-k}(\Sigma_g);
\mathbb{T}_\alpha, \mathbb{T}_\beta)$.

\begin{lemma} \label{lem2} Let $P \in \pi_1 (\Omega(\mathbb{T}_\alpha, \mathbb{T}_\beta), \mathbf{x})$ and $\iota(P)$ be the image of $P$ under the map defined in Lemma \ref{lem1}. Then we have the following equalities: 
\begin{eqnarray*}
	\text{Index}_{(\mathbb{T}_\alpha, \mathbb{T}_\beta)} (P)  &=& \text{Index}_{(\mathbb{L}_1, \mathbb{L}_2)} (\iota(P)) \\ 
        \text{Area}_{\sym^{g-k}(\Sigma_g)} (P) &=& \text{Area}_{\mathbb{M}} (\iota(P)) 
\end{eqnarray*}
\end{lemma}
\emph{Proof.} Recall that the map $\iota : \Omega(\mathbb{T}_\alpha, \mathbb{T}_\beta) \to \Omega(\mathbb{L}_1, \mathbb{L}_2) $ sends a path $\gamma \to (\gamma, x_2,x_3,\ldots, x_{2(g-k)})$ where $x_2,x_3,\ldots,x_{2(g-k)}$ are constant paths. Thus, if $P$ is a path of paths, written as $\gamma_s(t)$, the image of $P$ under $\iota$ is given by $(\gamma_s(t), x_2(s),x_3(s),\ldots, x_{2(g-k)}(s))$, hence, all but 
the first component has vanishing $t$-derivative. Therefore, the areas are the same as claimed since the area of the image of $P$ in $\f{M}$ has contribution only from the component mapping to $\sym^{g-k}(\Sigma_g)$. Similarly, the fact that Maslov-Viterbo index is the same follows form the fact that $\mathbb{T}_\alpha \cap \mathbb{T}_\beta = \mathbb{L}_1 \cap \mathbb{L}_2$ and crossing-forms used in calculation of the Maslov indices agree, this in turn follows from the fact that only the first component of $\iota(P)$ has non-vanishing $t$-derivative (see proof of Lemma 3.1.6 in \cite{WW} for a similar argument). \qed

In view of Lemma \ref{lem1} and Lemma \ref{lem2}, to conclude monotonicity for
the pair $(\mathbb{L}_1 , \mathbb{L}_2) $,  all we need to show is that index
and area are proportional for the pair of Lagrangians $(\mathbb{T}_\alpha,
\mathbb{T}_\beta)$ in $\sym^{g-k}(\Sigma_g)$. Until now, we have been working
with the symplectic forms $\omega$  on $\mathbb{M} = \prod_{b\in B}
\sym^{n_b} (F_{b})$ that Perutz constructed which have the cohomology class in
$\eta+\lambda \theta$ for sufficiently small $\lambda>0$ on each component. 
The reason for not using symplectic forms in the arguably simpler class $\eta$ is because Perutz's
construction that we appeal to in Theorem \ref{per} is only known to work for a
symplectic form in class $\eta+\lambda \theta$ for non-zero $\lambda$ (because
of the fact that the relative Hilbert scheme used in the construction is not
affine). 

Now, we claim that if we isotope $\alpha$ and $\beta$ curves in such way to
form an admissible diagram, then for a symplectic form $\omega_{\xi}$ (to be
made explicit below) in the class $\eta$ monotonicity is satisfied. By construction in Proposition 1.1 of \cite{handleslide} both $\omega$ and $\omega_\xi$ can be
arranged to agree outside a neighborhood of the diagonal and tame the same
regular almost complex structure $J_s=\sym^{g-k}(j_s)$ everywhere. This suffices
for the purpose of proving that Floer homology is well-defined. (It also follows that $QFH(Y,f)$ does not depend on the precise value of $\lambda$
	as long as it is sufficiently small). 

Let us recall the essential properties of $\omega_\xi$ and $\omega$ from
\cite{handleslide}. $\omega_\xi$ for $\xi$ an area form on $\Sigma_g$ is
characterized by the following equation. Let $P$ be a 2-dimensional region in
$\sym^{g-k}(\Sigma_g)$ and $D(P) \subset \Sigma_g$ be its projection to
$\Sigma_g$ defined by first taking the preimage of
$P$ by the map $\pi: \Sigma_g^{\times{(g-k)}} \to \sym^{g-k}(\Sigma_g)$ and
projecting to the first component (cf. \cite{OS} Definition 2.13). Then we have \[ \int_P \omega_\xi =
\frac{1}{(g-k)!} \int_{D(P)} \xi \]

Such symplectic forms $\omega_\xi$ were constructed by Perutz in $\cite{handleslide}$ and they represent the cohomology class $\pi_*([\xi^{\times (g-k)}]) = s \eta $ where $s = \frac{1}{(g-k)!} \int_{\Sigma_g} \xi$  and coincides with $\pi_*(\xi^{\times (g-k)})$ outside a neighborhood of the diagonal $\Delta$. We are free to scale $\xi$ if necessary in order to adjust $s$. 

The construction of the form $\omega$ in class $\eta+\lambda\theta$ in
$\cite{handleslide}$ is given by modifying $\omega_\xi$ near the diagonal.
Indeed $\omega = \omega_\xi - \epsilon \delta$ for $\epsilon >0$ small  where $\delta$ denotes a closed $(1,1)$ form representing
the diagonal class $[\Delta] = (4g-2k-2)\eta - 2\theta $ (\cite{macdonald})  and supported in a small neighborhood of $\Delta$. By taking $\epsilon$ sufficiently small, and adjusting $s$ by scaling $\xi$, we can obtain a symplectic form $\omega$ in class $\eta+ \lambda \theta$ for $\lambda>0$ sufficiently small. Then, it is easy to see that $\omega$ and $\omega_\xi$ tame the same regular complex structures $\sym^{g-k}(j_s)$ (since the modification is compactly supported and sufficiently small). 

Next, we adopt the definition of strong admissibility from $\cite{OS}$ to our setting. Note that as in \cite{OS} , $\mathbf{x}$ determines a spin$^c$ structure $\s
\in \mathcal{S}(Y|\Sigma_k)$. (In fact, again as in $\cite{OS}$ , one can see
that $\pi_0(\Omega(\mathbb{T}_\alpha, \mathbb{T}_\beta)) \simeq
\mathcal{S}(Y|\Sigma_k)$ ). Recall also that if $P \in
\pi_1(\Omega(\mathbb{T}_\alpha, \mathbb{T}_\beta), \mathbf{x})$ then we denote
the homology class $H(P) \in H_2(Y)$ determined by $P$ obtained as follows: $P$
can be projected onto the surface $\Sigma_g$ such that the projection bounds
various $\alpha_i$ and $\beta_j$, these boundary components then can be capped
off inside $Y$ by the corresponding cores of the handles attached to $\alpha_i$
and $\beta_j$ (see \cite{OS}).

\begin{definition} \label{defin} Let $(\Sigma_g, \alpha_1,\ldots, \alpha_{g-k},
	\beta_1,\ldots,\beta_{g-k})$ be a diagram obtained from a broken fibration on $Y$ and $k>1$, we say that the
	diagram is strongly admissible if there exists an area form $\xi$ on
	$\Sigma_g$ such that for any $P \in \pi_1(\Omega(\mathbb{T}_\alpha,
	\mathbb{T}_\beta), \mathbf{x})$ 
	the following holds: \[ (2-2k) \int_{D(P)} \xi = \langle c_1(\s), H(P) \rangle \int_{\Sigma_g} \xi \]       where $\s \in \mathcal{S}(Y|\Sigma_k)$ is the spin$^c$ structure determined by $\mathbf{x}$ and $[P] \in H_2(Y)$ is the homology class determined by $P$. 
\end{definition}

Note that splicing a sphere representing $\pi_2(\sym^{g-k}(\Sigma_g)) = \f{Z}$
, changes the index by \[ 2 \langle c_1(T\sym^{g-k})(\Sigma_g), [\Sigma_{g-k}]
\rangle = ((1-k)\eta-\theta) ([\Sigma_g]) = (2-2k)\] which is exactly equal to
$\langle c_1(\s), [\Sigma_g] \rangle$ for $\s \in \mathcal{S}(Y|\Sigma_k)$. This is why the normalization factor in our definition of strong admissibility differs form the one in $\cite{OS}$. 

Now, if we consider an arbitrary diagram $(\Sigma_g, \alpha_1,\ldots
\alpha_{g-k}, \beta_1,\ldots \beta_{g-k})$, by winding (say $\alpha$ curves) as
in \cite{OS} (cf Lemma \ref{admissible}), we can make the diagram strongly
admissible for $k>1$ (for $k=1$, we again need to use the Novikov ring
$\Lambda$ and keep track of the intersection with a basepoint $w$). The
following lemma shows that for $k>1$ strong admissibility implies monotonicity
for a symplectic form $\omega_\xi$ in the class $\eta$ (cf.
	\cite{lekiliPerutz}) and concludes the proof that $QFH(Y,f,\s)$ is
	well-defined for $\s \in \mathcal{S}(Y|\Sigma_k)$ and $g<2k$. 

\begin{proposition} For a strongly admissible diagram $(\Sigma_g,
	\alpha_1,\ldots,\alpha_{g-k},\beta_1,\ldots \beta_{g-k})$ ,
	the $\omega_\xi$-area and the index maps from
	$\pi_1(\Omega(\mathbb{T}_\alpha, \mathbb{T}_\beta),
	\mathbf{x}) \to \f{R}$ are proportional. \end{proposition}

	\emph{Proof.} If $P
	\in \pi_1(\Omega(\mathbb{T}_\alpha, \mathbb{T}_\beta), \mathbf{x}) $ ,
	it follows from the same argument as in Theorem 4.9 of \cite{OS} that
	$\text{Index}(P) = \langle c_1(\s), H(P) \rangle $, where $H(P) \in H_2(Y)$ is the homology class determined by $P$. Strong admissibility then implies that:
	\[ (2-2k) \int_{D(P)} \xi = \text{Index}(P) \int_{\Sigma_g} \xi \] \QED

\begin{corollary} When $k>1$, $QFH(Y,f,\s)$ is well-defined over $\f{Z}_2$ for $\s \in \mathcal{S}(Y|\Sigma_k)$ and $g<2k$.\end{corollary} \QED

\subsection{The isomorphism}

Recall that when defining $QFH'(Y,f,\Lambda)$ as a variant of Heegaard Floer homology
we have used Lipshitz's cylindrical reformulation, by setting up the theory in
$\Sigma_{\text{max}} \times [0,1] \times \f{R}$. This was convenient because of the
bubbling issues that may occur in the negatively monotone manifold
$\sym^{g-k}(\Sigma_{\text{max}})$. However, in the strongly negative case, when
$g<2k$, bubbling can be ruled out for a generic path $J_s$ of almost
complex structures on $\sym^{g-k}(\Sigma_{\text{max}})$ on the grounds that the moduli
space of bubbles in this case has negative virtual dimension. In fact, the proof of
the lemma below shows that in this case, one can also use a path of integrable complex
structures as in the case of the usual Heegaard Floer homology to make sense of this
group.  Thus, the Floer homology groups can be formulated as a Lagrangian intersection
theory in $\sym^{g-k}(\Sigma_{\text{max}})$.

\begin{lemma} \label{reformulation} Suppose that $Y$ admits a broken fibration with $
	g < 2k $. Then for $\s \in \mathcal{S}(Y|\Sigma_{\text{min}})$,
	\[ QFH'(Y,f,\s;\Lambda) \simeq HF(\sym^{g-k}(\Sigma_{\text{max}}); \alpha_1
	\times \ldots \times \alpha_{g-k}, \beta_{1} \times \ldots \times
	\beta_{g-k},\s ;\Lambda) \]

\end{lemma}

\emph{Proof.} We first argue that for a generic path of almost complex structures
$\{j_s\}$ on $\Sigma_{\text{max}}$ the induced integrable complex structures
$\sym^{g-k}(j_s)$ achieve transversality for the holomorphic disks mapping to
$\sym^{g-k}(\Sigma_{\text{max}})$ which contribute to the differential and furthermore
for these complex structures no bubbling can occur because of the strong negativity
assumption $g<2k$. The fact that these complex structures achieve transversality is
standard and follows exactly as in the case of the usual Heegaard Floer homology
set-up, see for example Proposition A.5 of \cite{lipshitz}. To avoid bubbling, we make
use of the Abel-Jacobi map: \[ \text{AJ} : \sym^{g-k}(\Sigma_{\text{max}}) \to
\text{Jac}(\Sigma_{\text{max}})\]  The assumption $g<2k$ ensures that the Abel-Jacobi
map is injective for $j$ chosen outside of a subset of complex codimension at least
$1$ (so that for a generic path $j_s$ it's injective for all $s$). A generic choice of
$j_s$ therefore ensures that there cannot be any non-constant holomorphic spheres mapping to
$\sym^{g-k}(\Sigma_{\text{max}})$ since $\pi_2(\text{Jac}(\Sigma_{\text{max}}))=0$. One can also rule out disk bubbles in the same way:
since the inclusions of $\alpha_1 \times\ldots\times \alpha_{g-k}$ and
$\beta_1\times\ldots\times\beta_{g-k}$ to $\sym^{g-k}(\Sigma_{\text{max}})$ are
injective at the level of fundamental groups and since the Abel-Jacobi map is
injective and induces an isomorphism on the first homology when $g<2k$, the image of
a holomorphic disc by the Abel-Jacobi map represents a trivial relative homology class, therefore it is trivial. Hence, there cannot be any non-constant holomorphic disk bubbles. 

Now, applying the reformulation of Lipshitz, as in \cite{lipshitz}, allows us to
translate the Lagrangian Floer homology in $\sym^{g-k}(\Sigma_{\text{max}})$
``tautologically'' to the cylindrical set-up in $\Sigma_\text{max} \times [0,1] \times \f{R}$
(see appendix A in \cite{lipshitz}).  \QED

Finally, we are ready to state our theorem that establishes the isomorphism between
quilted Floer homology groups arising from Lagrangian correspondences with Heegaard
Floer homology.  

\begin{theorem} \label{3mfldiso} Suppose that $Y$ admits a broken fibration with
	$g<2k$. Then for $\s \in
	\mathcal{S}(Y|\Sigma_{\text{min}})$, \[ HF^+(Y, f, \gamma_w, \s) \simeq
	QFH'(Y,f;\s, \Lambda) \simeq QFH(Y,f;\s, \Lambda) \] \end{theorem}

\emph{Proof.} The proof will be obtained by putting together the results obtained so
far together with the composition theorem for correspondences (\cite{LL},
\cite{katrin}) mentioned in the introduction to this section. This allows one to
compose Lagrangian correspondences and obtain isomorphic Floer homology groups.  More
precisely, Theorem \ref{isomorphism} and Lemma \ref{reformulation} give us that
$HF^+(Y,f, \gamma_w,\s) \simeq  QFH'(Y,f;\s,\Lambda) \simeq
HF(\sym^{g-k}(\Sigma_{\text{max}}); \alpha_1 \times \ldots \times \alpha_{g-k},
\beta_{1} \times \ldots \times \beta_{g-k};\Lambda)$. Now, Lemma \ref{lemham}
expresses the Lagrangians $\alpha_1 \times \ldots \times \alpha_{g-k}$ and $\beta_{1}
\times \ldots \times \beta_{g-k}$ as transverse and embedded compositions of the
Lagrangians $L_{\alpha_i}$ and $L_{\beta_j}$. Therefore, we are in a position to apply
the Wehrheim-Woodward's composition theorem which says that quilted Floer homology
groups associated with a cyclic set of Lagrangians correspondences is invariant under
transverse and embedded compositions of the Lagrangians. One important technicality
that arises in the proof of Wehrheim-Woodward is the possibility of ``figure-eight''
bubbles. To avoid those, Wehrheim-Woodward originally proved their theorem only in the
positively monotone case (and exact case), whereas we are in the strongly negative
setting and it is not clear that the approach of Wehrheim and Woodward can be
generalized to this situation. 

We have addressed this problem in another place (see \cite{LL}) where we gave a new proof of
Wehrheim-Woodward result which applies equally well in the strongly negative setting.

Therefore, we obtain an isomorphism between the Floer homology of the Lagrangians
$\alpha_1 \times \ldots \times \alpha_{g-k}$, $\beta_{1} \times \ldots \times
\beta_{g-k}$ and the quilted Floer homology of the Lagrangian correspondences
$L_{\alpha_1}, \ldots, L_{\alpha_{g-k}}$ and $L_{\beta_1}, \ldots, L_{\beta_{g-k}}$.
This completes the proof.\QED

\section{Discussion: 4--manifold invariants} \label{4manifolds}

We first recall the definition of broken Lefschetz fibrations on smooth $4$-manifolds.

\begin{definition}
A {\it broken fibration} on a closed $4$--manifold $X$ is a smooth map to a closed surface with singular set $A \cup B$, where $A$ is a finite set of singularities of Lefschetz type near which a local model in oriented charts is the complex map $(w, z ) \to w^2 + z^2$ , and $B$ is a $1$-dimensional submanifold along which the fibration is locally modelled by the real map $(t, x, y, z ) \to (t, x^2+y^2 -z^2 )$, $B$ corresponding to $t = 0$. 

\end{definition}

It was proven in \cite{lekili} that every closed oriented smooth $4$--manifold admits
an {\it equatorial} broken Lefschetz fibration to $S^2$ (see also \cite{AKkara} where
the authors give a new proof of this result using handlebody calculus). Equatorial here
means that the 1--dimensional part of the critical value set is a set of embedded
parallel circles on $S^2$.  Lagrangian matching invariants of a $4$--manifold as
defined by Perutz in \cite{LM1} are obtained by counting quilted holomorphic sections
of a broken fibration associated with the $4$--manifold. These invariants, which are
conjecturally equal to Seiberg-Witten invariants, have a TQFT-like structure where the
three manifold invariants are the quilted Floer homology groups that we have discussed
in this paper.  Similarly, Heegaard Floer homology is the three manifold part of a
TQFT-like structure, which underlies the construction of Ozsv\'ath-Szab\'o
$4$--manifold invariants \cite{OS4}. 

By cutting a broken fibration along a family of circles that are transverse to the
equatorial circles of critical values, one can obtain
a cobordism decomposition of the $4$--manifold, such that each cobordism is an
elementary cobordism, namely it is a cobordism obtained by either a one or two handle
attachment. Therefore, because of Theorem \ref{3mfldiso}, in order to equate the above
mentioned four-manifold invariants for the spin$^c$ structures which satisfy the
adjunction equality with respect to the minimal genus fibre of the broken fibration,
one needs to check only that the cobordism maps for one and two handle attachments in both
theories coincide. This will be in turn obtained by extending the techniques developed
in this paper to cobordism maps. We plan to investigate this latter claim in a sequel
to this paper. This will in particular prove that for the spin$^c$ structures
considered, the Lagrangian matching invariants are independent of the broken fibration
that is chosen on the $4$--manifold.

\end{document}